\documentclass[12pt]{article}
\usepackage{latexsym,amssymb,amsmath,enumerate,geometry,float,cite,tikz,setspace}
\newtheorem{definition}{Definition}
\newtheorem{theorem}[definition]{Theorem}

\newtheorem{lemma}[definition]{Lemma}

\newtheorem{claim}{Claim}

\begin{document}


\onehalfspace

\title{Dominating Sets inducing Large Components\\
in Maximal Outerplanar Graphs}

\author{Jos\'{e} D. Alvarado$^1$, Simone Dantas$^1$, and Dieter Rautenbach$^2$}

\date{}

\maketitle

\begin{center}
{\small 
$^1$ Instituto de Matem\'{a}tica e Estat\'{i}stica, Universidade Federal Fluminense, Niter\'{o}i, Brazil\\
\texttt{josealvarado.mat17@gmail.com, sdantas@im.uff.br}\\[3mm]
$^2$ Institute of Optimization and Operations Research, Ulm University, Ulm, Germany\\
\texttt{dieter.rautenbach@uni-ulm.de}
}
\end{center}

\begin{abstract}
For a maximal outerplanar graph $G$ of order $n$ at least $3$, 
Matheson and Tarjan showed that $G$ has domination number at most $n/3$.
Similarly, 
for a maximal outerplanar graph $G$ of order $n$ at least $5$, 
Dorfling, Hattingh, and Jonck showed, by a completely different approach, 
that $G$ has total domination number at most $2n/5$
unless $G$ is isomorphic to one of two exceptional graphs of order $12$.

We present a unified proof of a common generalization of these two results.
For every positive integer $k$, we specify a set ${\cal H}_k$ of graphs of order at least $4k+4$ and at most $4k^2-2k$ such that 
every maximal outerplanar graph $G$ of order $n$ at least $2k+1$ that does not belong to ${\cal H}_k$
has a dominating set $D$ of order at most $\lfloor\frac{kn}{2k+1}\rfloor$
such that every component of the subgraph $G[D]$ of $G$ induced by $D$ has order at least $k$.
\end{abstract}

{\small 

\medskip


\medskip


}

\pagebreak

\section{Introduction}\label{section1}

The two most prominent domination parameters \cite{hahesl,heye},
the domination number $\gamma(G)$
and 
the total domination number $\gamma_t(G)$
of a graph $G$,
have both been studied in detail for maximal outerplanar graphs \cite{aldara,cawa,cahemama,to}.
Two fundamental results in this context are as follows.

\begin{theorem}[Matheson and Tarjan \cite{mata}]\label{theoremmata}
If $G$ is a maximal outerplanar graph of order $n$ at least $3$, then $\gamma(G)\leq \left\lfloor\frac{n}{3}\right\rfloor$.
\end{theorem}

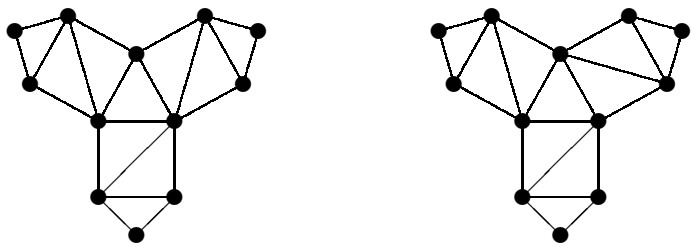
\begin{figure}[H]
\begin{center}
\unitlength 1mm 
\linethickness{0.4pt}
\ifx\plotpoint\undefined\newsavebox{\plotpoint}\fi 
\begin{picture}(33,30)(0,0)
\put(16,0){\circle*{2}}
\put(21,5){\circle*{2}}
\put(11,5){\circle*{2}}
\put(11,15){\circle*{2}}
\put(21,15){\circle*{2}}
\put(21,15){\line(0,-1){10}}
\put(21,5){\line(-1,0){10}}
\put(11,5){\line(1,-1){5}}
\put(16,0){\line(1,1){5}}
\put(21,15){\line(-1,-1){10}}
\put(11,5){\line(0,1){10}}
\put(21,15){\line(-1,0){10}}
\put(16,24){\circle*{2}}
\put(30,20){\circle*{2}}
\put(25,29){\circle*{2}}
\multiput(25,29)(-.060402685,-.033557047){149}{\line(-1,0){.060402685}}
\multiput(16,24)(.033557047,-.060402685){149}{\line(0,-1){.060402685}}
\multiput(21,15)(.060402685,.033557047){149}{\line(1,0){.060402685}}
\multiput(30,20)(-.033557047,.060402685){149}{\line(0,1){.060402685}}
\put(25,29){\line(0,1){0}}
\multiput(25,29)(-.033613445,-.117647059){119}{\line(0,-1){.117647059}}
\put(32,27){\circle*{2}}
\multiput(25,29)(.11666667,-.03333333){60}{\line(1,0){.11666667}}
\multiput(32,27)(-.03333333,-.11666667){60}{\line(0,-1){.11666667}}
\put(2,20){\circle*{2}}
\put(7,29){\circle*{2}}
\put(0,27){\circle*{2}}
\multiput(16,24)(-.060402685,.033557047){149}{\line(-1,0){.060402685}}
\multiput(7,29)(-.11666667,-.03333333){60}{\line(-1,0){.11666667}}
\multiput(0,27)(.03333333,-.11666667){60}{\line(0,-1){.11666667}}
\multiput(2,20)(.060402685,-.033557047){149}{\line(1,0){.060402685}}
\multiput(11,15)(.033557047,.060402685){149}{\line(0,1){.060402685}}
\multiput(7,29)(-.033557047,-.060402685){149}{\line(0,-1){.060402685}}
\multiput(7,29)(.033613445,-.117647059){119}{\line(0,-1){.117647059}}
\end{picture}\hspace{2cm}
\unitlength 1mm 
\linethickness{0.4pt}
\ifx\plotpoint\undefined\newsavebox{\plotpoint}\fi 
\begin{picture}(33,30)(0,0)
\put(16,0){\circle*{2}}
\put(21,5){\circle*{2}}
\put(11,5){\circle*{2}}
\put(11,15){\circle*{2}}
\put(21,15){\circle*{2}}
\put(21,15){\line(0,-1){10}}
\put(21,5){\line(-1,0){10}}
\put(11,5){\line(1,-1){5}}
\put(16,0){\line(1,1){5}}
\put(21,15){\line(-1,-1){10}}
\put(11,5){\line(0,1){10}}
\put(21,15){\line(-1,0){10}}
\put(16,24){\circle*{2}}
\put(30,20){\circle*{2}}
\put(25,29){\circle*{2}}
\multiput(25,29)(-.060402685,-.033557047){149}{\line(-1,0){.060402685}}
\multiput(16,24)(.033557047,-.060402685){149}{\line(0,-1){.060402685}}
\multiput(21,15)(.060402685,.033557047){149}{\line(1,0){.060402685}}
\multiput(30,20)(-.033557047,.060402685){149}{\line(0,1){.060402685}}
\put(25,29){\line(0,1){0}}
\put(32,27){\circle*{2}}
\multiput(25,29)(.11666667,-.03333333){60}{\line(1,0){.11666667}}
\multiput(32,27)(-.03333333,-.11666667){60}{\line(0,-1){.11666667}}
\put(2,20){\circle*{2}}
\put(7,29){\circle*{2}}
\put(0,27){\circle*{2}}
\multiput(16,24)(-.060402685,.033557047){149}{\line(-1,0){.060402685}}
\multiput(7,29)(-.11666667,-.03333333){60}{\line(-1,0){.11666667}}
\multiput(0,27)(.03333333,-.11666667){60}{\line(0,-1){.11666667}}
\multiput(2,20)(.060402685,-.033557047){149}{\line(1,0){.060402685}}
\multiput(11,15)(.033557047,.060402685){149}{\line(0,1){.060402685}}
\multiput(7,29)(-.033557047,-.060402685){149}{\line(0,-1){.060402685}}
\multiput(7,29)(.033613445,-.117647059){119}{\line(0,-1){.117647059}}
\multiput(16,24)(.117647059,-.033613445){119}{\line(1,0){.117647059}}
\end{picture}
\end{center}
\caption{The two exceptional graphs $H_1$ and $H_2$ for Theorem \ref{theoremdohajo}.}\label{fig1}
\end{figure}

\begin{theorem}[Dorfling, Hattingh, and Jonck \cite{dohajo}]\label{theoremdohajo}
If $G$ is a maximal outerplanar graph of order $n$ at least $5$ 
that is not isomorphic to one of the two graphs $H_1$ and $H_2$ in Figure \ref{fig1}, 
then $\gamma(G)\leq \left\lfloor\frac{2n}{5}\right\rfloor$.
\end{theorem}

\noindent The proofs of these two results in \cite{mata,dohajo} are quite different.
While Theorem \ref{theoremmata} follows from an elegant labeling argument, 
the proof of Theorem \ref{theoremdohajo} relied on a detailed case analysis;
one reason for this difference probably being the existence of the two exceptional graphs.

Our goal in the present paper is a unified proof of a common generalization of these two results.

\medskip

\noindent For some positive integer $k$ and a graph $G$,
a set $D$ of vertices of $G$ is a {\it $k$-component dominating set} of $G$
if every vertex in $V(G)\setminus D$ has a neighbor in $D$, and
every component of the subgraph $G[D]$ of $G$ induced by $D$ has order at least $k$.
The minimum cardinality of a $k$-component dominating set of $G$ 
is the {\it $k$-component domination number} $\gamma_k(G)$ of $G$.

Note that a graph has a $k$-component dominating set if and only if each of its components has order at least $k$.
Clearly, 
$\gamma_1(G)$ coincides with the domination number of $G$, 
and 
$\gamma_2(G)$ coincides with the total domination number of $G$, respectively.
The notation ``$\gamma_k(G)$'' has already been used 
to denote various other domination parameters. 
We chose this notation for its simplicity, 
and because there is little danger of confusion within the context of this paper.

For every positive integer $k$, we will specify a set ${\cal H}_k$ of graphs 
each of order at least $4k+4$ and at most $4k^2-2k$
such that our main result reads as follows.

\begin{theorem}\label{theorem1}
If $k$ and $n$ are positive integers with $n\geq 2k+1$, and $G$ is a maximal outerplanar graph of order $n$, then 
$$\gamma_k(G)\leq
\left\{
\begin{array}{ll}
\left\lceil\frac{kn}{2k+1}\right\rceil, & \mbox{if }G\in {\cal H}_k\\[3mm]
\left\lfloor\frac{kn}{2k+1}\right\rfloor, & \mbox{ otherwise.}
\end{array}
\right.$$
\end{theorem}
As we show below, the bound in Theorem \ref{theorem1} is actually tight for all values of $k$ and $n$ with $n\geq 2k+1$.
For $k=1$, we have $4k+4>4k^2-2k$, which implies that ${\cal H}_1$ is necessarily empty,
that is, Theorem \ref{theorem1} implies Theorem \ref{theoremmata}.
Furthermore, we will see that ${\cal H}_2=\{ H_1,H_2\}$,
that is, Theorems \ref{theorem1} implies Theorem \ref{theoremdohajo}.

The rest of the paper is devoted to the proof of Theorem \ref{theorem1}.

\section{Results}\label{section2}

For every maximal outerplanar graph, 
we will tacitly assume that it is embedded in the plane
in such a way that all its vertices are on the boundary of the unbounded face.
This implies that every bounded face is bounded by a triangle.
Furthermore, we assume that subgraphs 
inherit their embeddings in the natural way.

Let $G$ be a maximal outerplanar graph of order at least $3$.
The boundary of the unbounded face of $G$ is a Hamiltonian cycle $C(G)$ of $G$.
A chord of $G$ is an edge of $G$ that does not belong to $C(G)$.
Adding a chord $xy$ of $G$ to $C(G)$ results in a graph that has exactly two cycles $C_1$ and $C_2$ that are distinct from $C(G)$.
Furthermore, $C_1$ and $C_2$ are the boundaries of two maximal outerplanar subgraphs of $G$
whose union is $G$ and whose intersection is the edge $xy$.
We will refer to these two graphs as the {\it subgraphs of $G$ generated by $xy$}.
We refer to the edges of some graph $G$ as {\it $G$-edges}.

For positive integers $s$ and $t$, let $[s,t]$ be the set of positive integers at least $s$ and at most $t$, and let $[t]=[1,t]$.

\medskip

\noindent For positive integers $k$ and $n$ with $n\geq \max\{ 3,k\}$, let
$$\gamma_k(n)=\max\{ \gamma_k(G):\mbox{ $G$ is a maximal outerplanar graph of order $n$}\}.$$

\begin{lemma}\label{lemma1}
If $k$, $k'$, and $n$ are positive integers with $n\geq \max\{ 3,k\}$ and $k\geq k'$, then
\begin{enumerate}[(i)]
\item $\gamma_k(n)\leq \gamma_k(n+1)$, and
\item $\gamma_{k'}(n)\leq \gamma_k(n)$.
\end{enumerate}
\end{lemma}
{\it Proof:} (i) Let $G$ be a maximal outerplanar graph of order $n$ such that $\gamma_k(G)=\gamma_k(n)$. 
For some $C(G)$-edge $uv$ of $G$, 
let $G'$ arise from $G$ by adding a new vertex $x$ and the new edges $ux$ and $xv$.
Clearly, $G'$ is a maximal outerplanar graph of order $n+1$.
Let $D'$ be a minimum $k$-component dominating set of $G'$.
If either $x\not\in D'$ or $x\in D'$ and the component of $G'[D']$ that contains $x$ has order at least $k+1$,
then $D'\setminus \{ x\}$ is a $k$-component dominating set of $G$,
which implies
$\gamma_k(n)=\gamma_k(G)\leq\gamma_k(G')\leq \gamma_k(n+1)$.
Hence, we may assume that $x\in D'$, and 
that the component of $G'[D']$ that contains $x$ has order exactly $k$.
Since $n\geq k$, there is a vertex $y$ in $V(G')\setminus D'$ 
that has a neighbor in the component of $G'[D']$ that contains $x$.
The set $D=(D'\setminus \{ x\})\cup \{ y\}$ is a $k$-component dominating set of $G$,
which implies $\gamma_k(n)\leq \gamma_k(n+1)$ as above.

\medskip

\noindent (ii) This follows immediately from the trivial fact 
that every $k$-component dominating set is a $k'$-component dominating set. $\Box$

\begin{lemma}\label{lemma2}
If $k$ is a positive integer, then $\gamma_k(2k+3)=k$.
\end{lemma}
{\it Proof:} Since $\gamma_k(2k+3)\geq k$ follows immediately from the definition,
it remains to show $\gamma_k(2k+3)\leq k$,
which we prove by induction on $k$. 
Since every maximal outerplanar graph of order $5$ has a universal vertex, we obtain $\gamma_1(5)\leq 1$.
Now, let $k\geq 2$.
Let $G$ be a maximal outerplanar graph of order $2k+3$.
Let $x$ be a vertex of degree $2$ in $G$.
The neighbors of $x$ in $G$, say $u$ and $v$, are adjacent.
Let $G'$ arise from $G$ by removing $x$ and contracting the edge $uv$ to a new vertex $u^*$.
The order of $G'$ is $2(k-1)+3$.
By induction, $G'$ has a $(k-1)$-component dominating set $D'$ of order $k-1$.
If $u^*\in D'$, then let $D=(D'\setminus \{ u^*\})\cup \{ u,v\}$.
If $u^*\not\in D'$, then let $D$ arise from $D'$ by adding one vertex from $\{ u,v\}$ that has a neighbor in $D'$.
Note that $D$ is well defined because $D'$ is a dominating set of $G'$.
In both cases, $D$ is a $k$-component dominating set of $G$ of order $k$,
which completes the proof. $\Box$

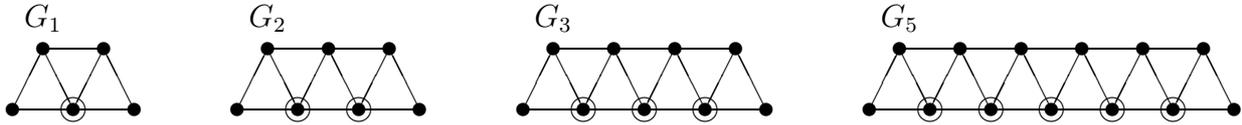
\begin{figure}[H]
\begin{center}
\unitlength 0.8mm 
\linethickness{0.4pt}
\ifx\plotpoint\undefined\newsavebox{\plotpoint}\fi 
\begin{picture}(21,15)(0,0)
\put(0,0){\circle*{2}}
\put(10,0){\circle*{2}}
\put(20,0){\circle*{2}}
\put(5,10){\circle*{2}}
\put(15,10){\circle*{2}}
\put(0,0){\line(1,2){5}}
\put(5,10){\line(1,-2){5}}
\put(10,0){\line(-1,0){10}}
\put(5,10){\line(1,0){10}}
\put(15,10){\line(-1,-2){5}}
\put(10,0){\line(1,0){10}}
\put(20,0){\line(-1,2){5}}
\put(5,15){\makebox(0,0)[cc]{$G_1$}}
\put(10,0){\circle{4}}
\end{picture}\hspace{1cm}
\linethickness{0.4pt}
\ifx\plotpoint\undefined\newsavebox{\plotpoint}\fi 
\begin{picture}(31,15)(0,0)
\put(0,0){\circle*{2}}
\put(10,0){\circle*{2}}
\put(20,0){\circle*{2}}
\put(30,0){\circle*{2}}
\put(5,10){\circle*{2}}
\put(15,10){\circle*{2}}
\put(25,10){\circle*{2}}
\put(0,0){\line(1,2){5}}
\put(5,10){\line(1,-2){5}}
\put(10,0){\line(-1,0){10}}
\put(5,10){\line(1,0){10}}
\put(15,10){\line(1,0){10}}
\put(15,10){\line(-1,-2){5}}
\put(25,10){\line(-1,-2){5}}
\put(10,0){\line(1,0){10}}
\put(20,0){\line(1,0){10}}
\put(20,0){\line(-1,2){5}}
\put(30,0){\line(-1,2){5}}
\put(5,15){\makebox(0,0)[cc]{$G_2$}}
\put(10,0){\circle{4}}
\put(20,0){\circle{4}}
\end{picture}\hspace{1cm}
\linethickness{0.4pt}
\ifx\plotpoint\undefined\newsavebox{\plotpoint}\fi 
\begin{picture}(41,15)(0,0)
\put(0,0){\circle*{2}}
\put(10,0){\circle*{2}}
\put(20,0){\circle*{2}}
\put(30,0){\circle*{2}}
\put(40,0){\circle*{2}}
\put(5,10){\circle*{2}}
\put(15,10){\circle*{2}}
\put(25,10){\circle*{2}}
\put(35,10){\circle*{2}}
\put(0,0){\line(1,2){5}}
\put(5,10){\line(1,-2){5}}
\put(10,0){\line(-1,0){10}}
\put(5,10){\line(1,0){10}}
\put(15,10){\line(1,0){10}}
\put(25,10){\line(1,0){10}}
\put(15,10){\line(-1,-2){5}}
\put(25,10){\line(-1,-2){5}}
\put(35,10){\line(-1,-2){5}}
\put(10,0){\line(1,0){10}}
\put(20,0){\line(1,0){10}}
\put(30,0){\line(1,0){10}}
\put(20,0){\line(-1,2){5}}
\put(30,0){\line(-1,2){5}}
\put(40,0){\line(-1,2){5}}
\put(5,15){\makebox(0,0)[cc]{$G_3$}}
\put(10,0){\circle{4}}
\put(20,0){\circle{4}}
\put(30,0){\circle{4}}
\end{picture}\hspace{1cm}
\linethickness{0.4pt}
\ifx\plotpoint\undefined\newsavebox{\plotpoint}\fi 
\begin{picture}(61,15)(0,0)
\put(0,0){\circle*{2}}
\put(10,0){\circle*{2}}
\put(20,0){\circle*{2}}
\put(30,0){\circle*{2}}
\put(40,0){\circle*{2}}
\put(50,0){\circle*{2}}
\put(60,0){\circle*{2}}
\put(5,10){\circle*{2}}
\put(15,10){\circle*{2}}
\put(25,10){\circle*{2}}
\put(35,10){\circle*{2}}
\put(45,10){\circle*{2}}
\put(55,10){\circle*{2}}
\put(0,0){\line(1,2){5}}
\put(5,10){\line(1,-2){5}}
\put(10,0){\line(-1,0){10}}
\put(5,10){\line(1,0){10}}
\put(15,10){\line(1,0){10}}
\put(25,10){\line(1,0){10}}
\put(35,10){\line(1,0){10}}
\put(45,10){\line(1,0){10}}
\put(15,10){\line(-1,-2){5}}
\put(25,10){\line(-1,-2){5}}
\put(35,10){\line(-1,-2){5}}
\put(45,10){\line(-1,-2){5}}
\put(55,10){\line(-1,-2){5}}
\put(10,0){\line(1,0){10}}
\put(20,0){\line(1,0){10}}
\put(30,0){\line(1,0){10}}
\put(40,0){\line(1,0){10}}
\put(50,0){\line(1,0){10}}
\put(20,0){\line(-1,2){5}}
\put(30,0){\line(-1,2){5}}
\put(40,0){\line(-1,2){5}}
\put(50,0){\line(-1,2){5}}
\put(60,0){\line(-1,2){5}}
\put(5,15){\makebox(0,0)[cc]{$G_5$}}
\put(10,0){\circle{4}}
\put(20,0){\circle{4}}
\put(30,0){\circle{4}}
\put(40,0){\circle{4}}
\put(50,0){\circle{4}}
\end{picture}
\end{center}
\caption{Some elements of a sequence $(G_k)_{k\in\mathbb{N}}$ 
of maximal outerplanar graphs of order $2k+3$.
Considering the two vertices of degree $2$ in the graph $G_{k+\ell}$ for positive integers $k$ and $\ell$ with $\ell\leq k$,
it follows easily that the encircled $k+\ell$ vertices form a minimum $k$-component dominating set of $G_{k+\ell}$.}\label{fig2}
\end{figure}

\begin{lemma}\label{lemma3}
If $k$ is a positive integer and $n\in [2k+1,4k+3]$, then $\gamma_k(n)=\lfloor\frac{kn}{2k+1}\rfloor$.
\end{lemma}
{\it Proof:} By definition and Lemmas \ref{lemma1} and \ref{lemma2}, we have 
$$k=\left\lfloor\frac{k(2k+1)}{2k+1}\right\rfloor
\leq \gamma_k(2k+1)
\leq \gamma_k(2k+2)
\leq \gamma_k(2k+3)
=\left\lfloor\frac{k(2k+3)}{2k+1}\right\rfloor=k,$$
which implies the statement for $2k+1\leq n\leq 2k+3$.

Now, let $n\in [2k+4,4k+3]$.
If $n$ is odd, say $n=2(k+\ell)+3$ for some positive integer $\ell$ with $\ell\leq k$,
then the graph $G_{k+\ell}$ illustrated in Figure \ref{fig2} easily implies
$\gamma_k(n)\geq \gamma_k(G_{k+\ell})=k+\ell$.
If $n$ is even, say $n=2(k+\ell)+2$ for some positive integer $\ell$ with $\ell\leq k$,
then the graph $G'_{k+\ell}$ that arises from $G_{k+\ell}$ by removing one vertex of degree $2$
easily implies $\gamma_k(n)\geq \gamma_k(G'_{k+\ell})=k+\ell$.
Since 
$$k+\ell\leq \frac{k(2(k+\ell)+2)}{2k+1}<\frac{k(2(k+\ell)+3)}{2k+1}<k+\ell+1$$
for $\ell\in [k]$,
this implies $\gamma_k(n)\geq \lfloor\frac{kn}{2k+1}\rfloor$.

For $\ell\in [k]$,
Lemmas \ref{lemma1} and \ref{lemma2} imply
$$\gamma_k(2(k+\ell)+2)\leq \gamma_k(2(k+\ell)+3)\leq \gamma_{k+\ell}(2(k+\ell)+3)=k+\ell,$$
which implies $\gamma_k(n)\leq \lfloor\frac{kn}{2k+1}\rfloor$, and completes the proof.
$\Box$

\medskip

\noindent It is an immediate consequence of Lemmas \ref{lemma2} and \ref{lemma3}
that for a positive integer $k$ and a non-negative integer $\ell$ with $\ell\leq k$, we have
\begin{eqnarray}\label{e1}
\gamma_k(2(k+\ell)+2)=\gamma_k(2(k+\ell)+3)=k+\ell.
\end{eqnarray}

\begin{lemma}\label{lemma4}
Let $k$ and $n$ be positive integers, and let $G$ be a maximal outerplanar graph of order $n$.
Let $u$ be a vertex of $G$, and let $xy$ be a $C(G)$-edge.
\begin{enumerate}[(i)]
\item If $n=2k+1$, then $G$ has a $k$-component dominating set $D$ of order $k$ that contains $u$.  
\item If $n=2k+2$, then $G$ has a $k$-component dominating set $D$ of order $k$ that intersects $xy$.
\item If $n=2k+1$, and $x$ and $y$ both have degree at least $3$ in $G$,
then $G$ has a $k$-component dominating set $D$ of order $k$ that contains $x$ and $y$.
\item If $n=2k+2$, and $x$ has degree at least $3$ in $G$,
then $G$ has a $k$-component dominating set $D$ of order $k$ that contains $x$.
\end{enumerate} 
\end{lemma}
{\it Proof:} Since the statements are trivial for $k=1$, we consider $k\geq 2$.

\medskip

\noindent (i) Since $2k+1=2(k-1)+3$, Lemma \ref{lemma2} implies that $G$ has $(k-1)$-component dominating set $D'$ of order $k-1$.
If $u\not\in D'$, then let $D=D'\cup \{ u\}$.
If $u\in D'$, then let $D=D'\cup \{ v\}$ where $v\in V(G)\setminus D'$.
In both cases, $D$ has the desired properties.

\medskip

\noindent (ii) Let $G'$ arise from $G$ by contracting the edge $xy$ to a new vertex $u^*$.
Since $G'$ has order $2k+1=2(k-1)+3$,
Lemma \ref{lemma2} implies that $G'$ has $(k-1)$-component dominating set $D'$ of order $k-1$.
If $u^*\in D'$, then let $D=(D'\setminus \{ u^*\})\cup \{ x,y\}$.
If $u^*\not\in D'$, then let $D$ arise from $D'$ by adding one vertex from $\{ x,y\}$ that has a neighbor in $D'$.
In both cases, $D$ has the desired properties.

\medskip

\noindent (iii) Let $z$ be a vertex of $G$ such that $xyz$ is a triangle in $G$.
Let $G_x$ be the subgraph of $G$ generated by the chord $yz$ that does not contain $x$,
and 
let $G_y$ be the subgraph of $G$ generated by the chord $xz$ that does not contain $y$.
Let $G_x$ and $G_y$ have orders $\ell_x+1$ and $\ell_y+1$, respectively.
Note that $n=\ell_x+\ell_y+1$, 
which implies that $\ell_x$ and $\ell_y$ have the same parity modulo $2$.

If $\ell_x$ and $\ell_y$ are both even, then (i) implies that 
$G_x$ has a $\ell_x/2$-component dominating set $D_x$ of order $\ell_x/2$ that contains $y$,
and  
$G_y$ has a $\ell_y/2$-component dominating set $D_y$ of order $\ell_y/2$ that contains $x$.
Since $k=\frac{n-1}{2}=\ell_x/2+\ell_y/2$,
possibly adding one vertex to the set $D_x\cup D_y$ yields a set with the desired properties.

If $\ell_x$ and $\ell_y$ are both odd, then (ii) implies that 
$G_x$ has a $(\ell_x-1)/2$-component dominating set $D_x$ of order $(\ell_x-1)/2$ that intersects $yz$,
and  
$G_y$ has a $(\ell_y-1)/2$-component dominating set $D_y$ of order $(\ell_y-1)/2$ that intersects $xz$.
The set $D=D_x\cup D_y\cup \{ x,y\}$ is a dominating set of $G$
such that $G[D]$ is connected and $|D|\leq |D_x|+|D_y|+1=k$.
Possibly adding further vertices to $D$ yields a set with the desired properties.

\medskip

\noindent (iv) If $y$ has degree $2$ in $G$, then $G'=G-y$ is a maximal outerplanar graph of order $2k+1$.
By (i), $G'$ has a $k$-component dominating set $D'$ of order $k$ that contains $x$.
Clearly, $D'$ is also a $k$-component dominating set of $G$.
Hence, we may assume that $y$ has degree at least $3$ in $G$.
Let $z$, $G_x$, $G_y$, $\ell_x$, and $\ell_y$ be as in (iii).
Since $n$ is even, $\ell_x$ and $\ell_y$ have different parities modulo $2$.

If $\ell_x$ is odd and $\ell_y$ is even, 
then (i) implies that 
$G_y$ has a $\ell_y/2$-component dominating set $D_y$ of order $\ell_y/2$ that contains $x$,
and (ii) implies that
$G_x$ has a $(\ell_x-1)/2$-component dominating set $D_x$ of order $(\ell_x-1)/2$ that intersects $yz$.
Since $(\ell_x+\ell_y-1)/2=(n-2)/2=k$,
possibly adding one further vertex to $D_x\cup D_y$ yields a set $D$ with the desired properties.

If $\ell_x$ is even and $\ell_y$ is odd, 
then (ii) implies that 
$G_y$ has a $(\ell_y-1)/2$-component dominating set $D_y$ of order $(\ell_y-1)/2$ that intersects $xz$,
and (i) implies that
$G_x$ has a $\ell_x/2$-component dominating set $D_x$ of order $\ell_x/2$ that contains $z$.
If $x\in D_y$, then possibly adding one vertex to $D_x\cup D_y$ yields a set $D$ with the desired properties.
If $x\not\in D_y$, then $z\in D_x\cap D_y$, and $D=D_x\cup D_y\cup\{ x\}$ has the desired properties.
$\Box$

\medskip

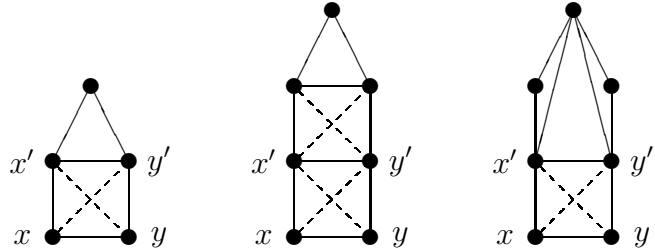
\begin{figure}[H]
\begin{center}
\unitlength 1mm 
\linethickness{0.4pt}
\ifx\plotpoint\undefined\newsavebox{\plotpoint}\fi 
\begin{picture}(19,21)(0,0)
\put(5,0){\circle*{2}}
\put(5,10){\circle*{2}}
\put(15,0){\circle*{2}}
\put(15,10){\circle*{2}}
\put(15,0){\line(-1,0){10}}
\put(15,10){\line(-1,0){10}}
\put(10,20){\circle*{2}}
\put(15,0){\line(0,1){10}}
\put(15,10){\line(-1,2){5}}
\put(10,20){\line(-1,-2){5}}
\put(5,10){\line(0,-1){10}}
\multiput(14.93,-.07)(-.0328947,.0328947){19}{\line(0,1){.0328947}}
\multiput(13.68,1.18)(-.0328947,.0328947){19}{\line(0,1){.0328947}}
\multiput(12.43,2.43)(-.0328947,.0328947){19}{\line(0,1){.0328947}}
\multiput(11.18,3.68)(-.0328947,.0328947){19}{\line(-1,0){.0328947}}
\multiput(9.93,4.93)(-.0328947,.0328947){19}{\line(0,1){.0328947}}
\multiput(8.68,6.18)(-.0328947,.0328947){19}{\line(0,1){.0328947}}
\multiput(7.43,7.43)(-.0328947,.0328947){19}{\line(0,1){.0328947}}
\multiput(6.18,8.68)(-.0328947,.0328947){19}{\line(0,1){.0328947}}
\multiput(14.93,9.93)(-.0328947,-.0328947){19}{\line(0,-1){.0328947}}
\multiput(13.68,8.68)(-.0328947,-.0328947){19}{\line(0,-1){.0328947}}
\multiput(12.43,7.43)(-.0328947,-.0328947){19}{\line(0,-1){.0328947}}
\multiput(11.18,6.18)(-.0328947,-.0328947){19}{\line(0,-1){.0328947}}
\multiput(9.93,4.93)(-.0328947,-.0328947){19}{\line(0,-1){.0328947}}
\multiput(8.68,3.68)(-.0328947,-.0328947){19}{\line(0,-1){.0328947}}
\multiput(7.43,2.43)(-.0328947,-.0328947){19}{\line(0,-1){.0328947}}
\multiput(6.18,1.18)(-.0328947,-.0328947){19}{\line(0,-1){.0328947}}
\put(4.93,-.07){\line(0,-1){.5}}
\put(1,0){\makebox(0,0)[cc]{$x$}}
\put(19,0){\makebox(0,0)[cc]{$y$}}
\put(19,10){\makebox(0,0)[cc]{$y'$}}
\put(1,10){\makebox(0,0)[cc]{$x'$}}
\end{picture}\hspace{1cm}
\linethickness{0.4pt}
\ifx\plotpoint\undefined\newsavebox{\plotpoint}\fi 
\begin{picture}(19,31)(0,0)
\put(5,0){\circle*{2}}
\put(5,10){\circle*{2}}
\put(5,20){\circle*{2}}
\put(15,0){\circle*{2}}
\put(15,10){\circle*{2}}
\put(15,20){\circle*{2}}
\put(15,0){\line(-1,0){10}}
\put(15,10){\line(-1,0){10}}
\put(15,20){\line(-1,0){10}}
\put(10,30){\circle*{2}}
\put(15,0){\line(0,1){10}}
\put(15,10){\line(0,1){10}}
\put(15,20){\line(-1,2){5}}
\put(10,30){\line(-1,-2){5}}
\put(5,10){\line(0,-1){10}}
\put(5,20){\line(0,-1){10}}
\multiput(14.93,-.07)(-.0328947,.0328947){19}{\line(0,1){.0328947}}
\multiput(13.68,1.18)(-.0328947,.0328947){19}{\line(0,1){.0328947}}
\multiput(12.43,2.43)(-.0328947,.0328947){19}{\line(0,1){.0328947}}
\multiput(11.18,3.68)(-.0328947,.0328947){19}{\line(-1,0){.0328947}}
\multiput(9.93,4.93)(-.0328947,.0328947){19}{\line(0,1){.0328947}}
\multiput(8.68,6.18)(-.0328947,.0328947){19}{\line(0,1){.0328947}}
\multiput(7.43,7.43)(-.0328947,.0328947){19}{\line(0,1){.0328947}}
\multiput(6.18,8.68)(-.0328947,.0328947){19}{\line(0,1){.0328947}}
\multiput(14.93,9.93)(-.0328947,.0328947){19}{\line(0,1){.0328947}}
\multiput(13.68,11.18)(-.0328947,.0328947){19}{\line(0,1){.0328947}}
\multiput(12.43,12.43)(-.0328947,.0328947){19}{\line(0,1){.0328947}}
\multiput(11.18,13.68)(-.0328947,.0328947){19}{\line(0,1){.0328947}}
\multiput(9.93,14.93)(-.0328947,.0328947){19}{\line(0,1){.0328947}}
\multiput(8.68,16.18)(-.0328947,.0328947){19}{\line(0,1){.0328947}}
\multiput(7.43,17.43)(-.0328947,.0328947){19}{\line(0,1){.0328947}}
\multiput(6.18,18.68)(-.0328947,.0328947){19}{\line(0,1){.0328947}}
\multiput(14.93,9.93)(-.0328947,-.0328947){19}{\line(0,-1){.0328947}}
\multiput(13.68,8.68)(-.0328947,-.0328947){19}{\line(0,-1){.0328947}}
\multiput(12.43,7.43)(-.0328947,-.0328947){19}{\line(0,-1){.0328947}}
\multiput(11.18,6.18)(-.0328947,-.0328947){19}{\line(0,-1){.0328947}}
\multiput(9.93,4.93)(-.0328947,-.0328947){19}{\line(0,-1){.0328947}}
\multiput(8.68,3.68)(-.0328947,-.0328947){19}{\line(0,-1){.0328947}}
\multiput(7.43,2.43)(-.0328947,-.0328947){19}{\line(0,-1){.0328947}}
\multiput(6.18,1.18)(-.0328947,-.0328947){19}{\line(0,-1){.0328947}}
\multiput(14.93,19.93)(-.0328947,-.0328947){19}{\line(0,-1){.0328947}}
\multiput(13.68,18.68)(-.0328947,-.0328947){19}{\line(0,-1){.0328947}}
\multiput(12.43,17.43)(-.0328947,-.0328947){19}{\line(0,-1){.0328947}}
\multiput(11.18,16.18)(-.0328947,-.0328947){19}{\line(0,-1){.0328947}}
\multiput(9.93,14.93)(-.0328947,-.0328947){19}{\line(0,-1){.0328947}}
\multiput(8.68,13.68)(-.0328947,-.0328947){19}{\line(0,-1){.0328947}}
\multiput(7.43,12.43)(-.0328947,-.0328947){19}{\line(0,-1){.0328947}}
\multiput(6.18,11.18)(-.0328947,-.0328947){19}{\line(0,-1){.0328947}}
\put(4.93,-.07){\line(0,-1){.5}}
\put(1,0){\makebox(0,0)[cc]{$x$}}
\put(19,0){\makebox(0,0)[cc]{$y$}}
\put(19,10){\makebox(0,0)[cc]{$y'$}}
\put(1,10){\makebox(0,0)[cc]{$x'$}}
\end{picture}\hspace{1cm}
\unitlength 1mm 
\linethickness{0.4pt}
\ifx\plotpoint\undefined\newsavebox{\plotpoint}\fi 
\begin{picture}(19,31)(0,0)
\put(5,0){\circle*{2}}
\put(5,10){\circle*{2}}
\put(5,20){\circle*{2}}
\put(15,0){\circle*{2}}
\put(15,10){\circle*{2}}
\put(15,20){\circle*{2}}
\put(15,0){\line(-1,0){10}}
\put(15,10){\line(-1,0){10}}
\put(10,30){\circle*{2}}
\put(15,0){\line(0,1){10}}
\put(15,10){\line(0,1){10}}
\put(15,20){\line(-1,2){5}}
\put(10,30){\line(-1,-2){5}}
\put(5,10){\line(0,-1){10}}
\put(5,20){\line(0,-1){10}}
\multiput(14.93,-.07)(-.0328947,.0328947){19}{\line(0,1){.0328947}}
\multiput(13.68,1.18)(-.0328947,.0328947){19}{\line(0,1){.0328947}}
\multiput(12.43,2.43)(-.0328947,.0328947){19}{\line(0,1){.0328947}}
\multiput(11.18,3.68)(-.0328947,.0328947){19}{\line(-1,0){.0328947}}
\multiput(9.93,4.93)(-.0328947,.0328947){19}{\line(0,1){.0328947}}
\multiput(8.68,6.18)(-.0328947,.0328947){19}{\line(0,1){.0328947}}
\multiput(7.43,7.43)(-.0328947,.0328947){19}{\line(0,1){.0328947}}
\multiput(6.18,8.68)(-.0328947,.0328947){19}{\line(0,1){.0328947}}
\multiput(14.93,9.93)(-.0328947,-.0328947){19}{\line(0,-1){.0328947}}
\multiput(13.68,8.68)(-.0328947,-.0328947){19}{\line(0,-1){.0328947}}
\multiput(12.43,7.43)(-.0328947,-.0328947){19}{\line(0,-1){.0328947}}
\multiput(11.18,6.18)(-.0328947,-.0328947){19}{\line(0,-1){.0328947}}
\multiput(9.93,4.93)(-.0328947,-.0328947){19}{\line(0,-1){.0328947}}
\multiput(8.68,3.68)(-.0328947,-.0328947){19}{\line(0,-1){.0328947}}
\multiput(7.43,2.43)(-.0328947,-.0328947){19}{\line(0,-1){.0328947}}
\multiput(6.18,1.18)(-.0328947,-.0328947){19}{\line(0,-1){.0328947}}
\put(4.93,-.07){\line(0,-1){.5}}
\put(1,0){\makebox(0,0)[cc]{$x$}}
\put(19,0){\makebox(0,0)[cc]{$y$}}
\put(19,10){\makebox(0,0)[cc]{$y'$}}
\put(1,10){\makebox(0,0)[cc]{$x'$}}
\put(15,10){\line(-1,4){5}}
\put(6,10){\line(0,1){0}}
\put(10,30){\line(-1,-4){5}}
\end{picture}
\end{center}
\caption{All graphs $G$ for which $(G,xy)$ belongs to ${\cal G}_4$ and ${\cal G}_6$.
From each pair of crossing dashed edges exactly one edge belongs to $G$.}\label{figq}
\end{figure}

\noindent For an even integer $\ell$ at least $4$, let ${\cal G}_\ell$ be the set of all pairs $(G,xy)$, where 
\begin{itemize}
\item $G$ is a maximal outerplanar graph of order $\ell+1$, 
\item $xy$ is a $C(G)$-edge such that $\{ d_G(x),d_G(y)\}=\{ 2,3\}$, and
\item if $N_{C(G)}(x)=\{ x',y\}$ and $N_{C(G)}(y)=\{ y',x\}$, 
then the maximal outerplanar graph $G^-=G-\{ x,y\}$
does not have a $(\ell/2-2)$-component dominating set of order $\ell/2-2$ that intersects $x'y'$.
\end{itemize}
See Figure \ref{figq} for an illustration.
In fact, generalizing the first two graphs in this figure in the obvious way
implies that ${\cal G}_\ell$ is non-empty for every even $\ell$ at least $4$.

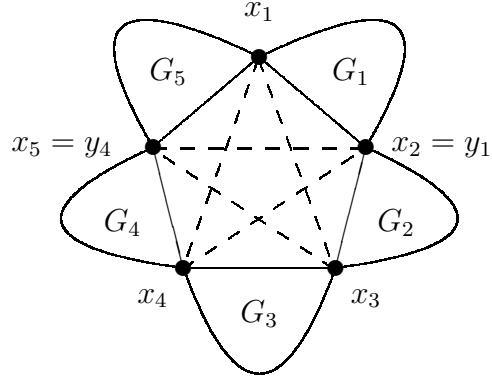
\begin{figure}[H]
\begin{center}
\unitlength 2mm 
\linethickness{0.4pt}
\ifx\plotpoint\undefined\newsavebox{\plotpoint}\fi 
\begin{picture}(40,28)(0,0)
\put(15,7){\circle*{1}}
\put(25,7){\circle*{1}}
\put(13,15){\circle*{1}}
\put(27,15){\circle*{1}}
\put(20,21){\circle*{1}}
\multiput(20,21)(-.039325843,-.033707865){178}{\line(-1,0){.039325843}}
\put(13,15){\line(1,-4){2}}
\put(15,7){\line(1,0){10}}
\put(25,7){\line(1,4){2}}
\multiput(27,15)(-.039325843,.033707865){178}{\line(-1,0){.039325843}}
\qbezier(20,21)(34.5,28)(27,15)
\qbezier(20,21)(5.5,28)(13,15)
\put(26,20){\makebox(0,0)[cc]{$G_1$}}
\put(29,10){\makebox(0,0)[cc]{$G_2$}}
\put(20,4){\makebox(0,0)[cc]{$G_3$}}
\put(11,10){\makebox(0,0)[cc]{$G_4$}}
\put(14,20){\makebox(0,0)[cc]{$G_5$}}
\multiput(19.93,20.93)(.03125,-.0875){10}{\line(0,-1){.0875}}
\multiput(20.555,19.18)(.03125,-.0875){10}{\line(0,-1){.0875}}
\multiput(21.18,17.43)(.03125,-.0875){10}{\line(0,-1){.0875}}
\multiput(21.805,15.68)(.03125,-.0875){10}{\line(0,-1){.0875}}
\multiput(22.43,13.93)(.03125,-.0875){10}{\line(0,-1){.0875}}
\multiput(23.055,12.18)(.03125,-.0875){10}{\line(0,-1){.0875}}
\multiput(23.68,10.43)(.03125,-.0875){10}{\line(0,-1){.0875}}
\multiput(24.305,8.68)(.03125,-.0875){10}{\line(0,-1){.0875}}
\multiput(24.93,6.93)(-.05,.0333333){15}{\line(-1,0){.05}}
\multiput(23.43,7.93)(-.05,.0333333){15}{\line(-1,0){.05}}
\multiput(21.93,8.93)(-.05,.0333333){15}{\line(-1,0){.05}}
\multiput(20.43,9.93)(-.05,.0333333){15}{\line(-1,0){.05}}
\multiput(18.93,10.93)(-.05,.0333333){15}{\line(-1,0){.05}}
\multiput(17.43,11.93)(-.05,.0333333){15}{\line(-1,0){.05}}
\multiput(15.93,12.93)(-.05,.0333333){15}{\line(-1,0){.05}}
\multiput(14.43,13.93)(-.05,.0333333){15}{\line(-1,0){.05}}
\put(26.93,14.93){\line(-1,0){.9333}}
\put(25.063,14.93){\line(-1,0){.9333}}
\put(23.196,14.93){\line(-1,0){.9333}}
\put(21.33,14.93){\line(-1,0){.9333}}
\put(19.463,14.93){\line(-1,0){.9333}}
\put(17.596,14.93){\line(-1,0){.9333}}
\put(15.73,14.93){\line(-1,0){.9333}}
\put(13.863,14.93){\line(-1,0){.9333}}
\multiput(19.93,20.93)(-.03125,-.0875){10}{\line(0,-1){.0875}}
\multiput(19.305,19.18)(-.03125,-.0875){10}{\line(0,-1){.0875}}
\multiput(18.68,17.43)(-.03125,-.0875){10}{\line(0,-1){.0875}}
\multiput(18.055,15.68)(-.03125,-.0875){10}{\line(0,-1){.0875}}
\multiput(17.43,13.93)(-.03125,-.0875){10}{\line(0,-1){.0875}}
\multiput(16.805,12.18)(-.03125,-.0875){10}{\line(0,-1){.0875}}
\multiput(16.18,10.43)(-.03125,-.0875){10}{\line(0,-1){.0875}}
\multiput(15.555,8.68)(-.03125,-.0875){10}{\line(0,-1){.0875}}
\multiput(14.93,6.93)(.05,.0333333){15}{\line(1,0){.05}}
\multiput(16.43,7.93)(.05,.0333333){15}{\line(1,0){.05}}
\multiput(17.93,8.93)(.05,.0333333){15}{\line(1,0){.05}}
\multiput(19.43,9.93)(.05,.0333333){15}{\line(1,0){.05}}
\multiput(20.93,10.93)(.05,.0333333){15}{\line(1,0){.05}}
\multiput(22.43,11.93)(.05,.0333333){15}{\line(1,0){.05}}
\multiput(23.93,12.93)(.05,.0333333){15}{\line(1,0){.05}}
\multiput(25.43,13.93)(.05,.0333333){15}{\line(1,0){.05}}
\put(32,15){\makebox(0,0)[cc]{$x_2=y_1$}}
\put(20,24){\makebox(0,0)[cc]{$x_1$}}
\put(27,5){\makebox(0,0)[cc]{$x_3$}}
\put(13,5){\makebox(0,0)[cc]{$x_4$}}
\put(7,15){\makebox(0,0)[cc]{$x_5=y_4$}}
\qbezier(25,7)(40,9)(27,15)
\qbezier(15,7)(0,9)(13,15)
\qbezier(15,7)(20,-7)(25,7)
\end{picture}
\end{center}
\caption{The embedding of a graph $G$ in ${\cal H}_k^2$. 
The dashed edges indicate a triangulation of the cycle $C_0(G):x_1x_2x_3x_5x_1$.}\label{figh}
\end{figure}

\noindent For positive integers $k$ and $p$ with $p\leq k-1$,
let ${\cal H}_k^p$ be the set of all graphs $G$ such that there are 
\begin{itemize}
\item $2p+1$ even integers $\ell_1,\ldots,\ell_{2p+1}$ with $4\leq \ell_i\leq 2k$ for $i\in [2p+1]$ and 
$$\ell_1+\cdots+\ell_{2p+1}\geq 4kp+2p+2,$$ 
as well as
\item $2p+1$ pairs 
$$(G_1,x_1y_1),\ldots,(G_{2p+1},x_{2p+1}y_{2p+1})$$ 
with $(G_i,x_iy_i)\in {\cal G}_{\ell_i}$ for $i\in [2p+1]$ such that
\end{itemize}
$G$ arises from the disjoint union of $G_1,\ldots,G_{2p+1}$ by 
\begin{itemize}
\item identifying the two vertices $y_i$ and $x_{i+1}$ for every $i\in [2p+1]$, 
where indices are identified modulo $2p+1$, and 
\item triangulating the cycle $C_0(G):x_1x_2\ldots x_{2p+1}x_1$.
\end{itemize}
The graphs in ${\cal H}_k^p$ have a natural embedding illustrated in Figure \ref{figh}.
In what follows, we always assume the graphs in ${\cal H}_k^p$ to be embedded in this way.

\noindent Let $${\cal H}_k=\bigcup_{p=1}^{k-1}{\cal H}_k^p.$$
Since for $k=1$, there is no positive integer $p$ with $p\leq k-1$, the set ${\cal H}_1$ is empty.
By definition, ${\cal H}_2={\cal H}_2^1$.
In view of the unique element of ${\cal G}_4$ illustrated in Figure \ref{figq},
the two graphs in Figure \ref{fig1} form the only elements of ${\cal H}_2$.

\begin{lemma}\label{lemma5}
Let $k$ and $p$ be positive integers with $p\leq k-1$.
Let $G$ be a graph in ${\cal H}_k^p$, and let $\ell_i$ and $(G_i,x_iy_i)$ for $i\in [2p+1]$ be as above.
\begin{enumerate}[(i)]
\item $G$ has even order $n$ at least $4kp+2p+2$ and at most $2k(2p+1)$.
\item $\ell_i+\ell_{i+1}\geq 2k+2p+2$ for every $i\in [2p+1]$.
\item $\gamma_k(G)=n/2-p=\left\lceil\frac{kn}{2k+1}\right\rceil$.
\item For every vertex $u$ in $V(G)\setminus V(C_0(G))$, the graph $G$ has a dominating set $D$ of order at most $n/2-(p+1)$ such that 
\begin{itemize}
\item $D$ has a partition into two disjoint sets $D_1$ and $D_2$,
\item $G[D_1]$ is a connected graph of order at least $\frac{\min\{ \ell_i:i\in [2p+1]\}}{2}-1$ that contains $u$,
\item every component of $G[D_2]$ has order at least $k$, and
\item $D$ contains no neighbor of $u$ on $C_0(G)$.
\end{itemize}
\end{enumerate}
\end{lemma}
{\it Proof:} (i) The lower bound on $n$ and its parity modulo $2$ are part of the definition of ${\cal H}_k^p$.
Since each $\ell_i$ is at most $2k$, the upper bound follows immediately.
 
\medskip

\noindent (ii) If there is some index $i\in [2p+1]$ with $\ell_i+\ell_{i+1}<2k+2k+2$, 
then, since, each of the remaining $2p-1$ values $\ell_j$ is at most $2k$,
we obtain $n<2k+2p+2+(2p-1)2k=4kp+2p+2$, which contradicts (i). 

\medskip

\noindent (iii) Note that $n/2-p=\left\lceil\frac{kn}{2k+1}\right\rceil$ is equivalent to 
$0=\left\lceil\frac{1}{2k+1}\Big((2k+1)p-n/2\Big)\right\rceil$,
which is equivalent to $-2k\leq (2k+1)p-n/2\leq 0$.
Therefore, the equality $n/2-p=\left\lceil\frac{kn}{2k+1}\right\rceil$
follows easily from (i). 

By Lemma \ref{lemma4}(i), for every $i\in [2p+1]$, 
the graph $G_i$ has 
a $\ell_i/2$-component dominating set $D_i^x$ of order $\ell_i/2$ that contains $x_i$
as well as 
a $\ell_i/2$-component dominating set $D_i^y$ of order $\ell_i/2$ that contains $y_i=x_{i+1}$.
By (ii), the set
$$D_1^y\cup D_2^x\cup D_3^y\cup D_4^x\cup \cdots\cup D_{2p-1}^y\cup D_{2p}^x\cup D_{2p+1}^y$$
is a $k$-component dominating set of $G$ of order at most $n/2-p$, which implies $\gamma_k(G)\leq n/2-p$.

It remains to show $\gamma_k(G)\geq n/2-p$.
Therefore, let $D$ be a $k$-component dominating set of minimum order of $G$ 
such that $|D\cap V(C_0(G))|$ is as large as possible.
For $i\in [2p+1]$, let $D_i=D\cap (V(G_i)\setminus \{ x_i,y_i\})$, 
$N_{C(G_i)}(x_i)=\{ x_i',y_i\}$, and 
$N_{C(G_i)}(y_i)=\{ x_i,y_i'\}$.

If there is some $i\in [2p+1]$ such that $|D_i|\leq \ell_i/2-2<k$, then, since $D$ is a $k$-component dominating set of $G$,
the set $D_i$ intersects $x_i'y_i'$. 
This implies that $D_i$ is a subset of some $(\ell_i/2-2)$-component dominating set 
of the graph $G_i^-=G_i-\{ x_i,y_i\}$ that is of order $\ell_i/2-2$ and intersects $x_i'y_i'$,
which implies the contradiction $(G_i,x_iy_i)\not\in {\cal G}_{\ell_i}$.
Hence $|D_i|\geq \ell_i/2-1$ for every $i\in [2p+1]$.

If there is some $i\in [2p+1]$ such that $x_i,x_{i+1}\not\in D$, then, since $D$ is a $k$-component dominating set of $G$,
the set $D_i$ has at least $k\geq \ell_i/2$ elements.
By Lemma \ref{lemma4} (i), 
the graph $G_i$ has a $k$-component dominating set $D_i'$ of order $k$ that contains $x_i$.
Now, $D'=(D\setminus D_i)\cup D_i'$ is a $k$-component dominating set of $G$ such that 
$|D'|\leq |D|$ and $|D'\cap V(C_0(G))|>|D\cap V(C_0(G))|$, 
which contradicts the choice of $D$.
Hence, for every $i\in [2p+1]$, we have $|D\cap \{ x_i,x_{i+1}\}|\geq 1$, which implies $|D\cap V(C_0(G))|\geq p+1$.

Altogether, we obtain 
$$|D|\geq \sum_{i=1}^{2p+1}(\ell_i/2-1)+(p+1)=n/2-(2p+1)+(p+1)=n/2-p,$$
which completes the proof of (iii).

\medskip

\noindent (iv) By symmetry, we may assume that $u\in V(G_{2p+1})$.
The graph $G^-_{2p+1}=G_{2p+1}-\{ x_{2p+1},y_{2p+1}\}$ has order $\ell_{2p+1}-1$.
By Lemma \ref{lemma4} (i), 
$G^-_{2p+1}$ has a $(\ell_{2p+1}/2-1)$-component dominating set $D^-_{2p+1}$ of order $\ell_{2p+1}/2-1$ that contains $u$.
Let $D_i^x$ and $D_i^y$ be as in (iii).
Let $D=D_1\cup D_2$, where 
$D_1=D^-_{2p+1}$ and 
$D_2=D_1^y\cup D_2^x\cup D_3^y\cup D_4^x\cup \cdots\cup D_{2p-1}^y\cup D_{2p}^x$.
Clearly, 
$D$ is a dominating set of $G$ of order at most $n/2-(p+1)$,
$D_1$ and $D_2$ are disjoint,
$G[D_1]$ is a connected graph of order $\ell_{2p+1}/2-1$ that contains $u$, and
every component of $G[D_2]$ has order at least $k$.
If $D$ contains a neighbor, say $v$, of $u$ on $C_0(G)$, then $v\in D_2$,
and it follows that the component of $G[D]$ that contains $u$ has order at least $k$,
which implies the contradiction that $D$ is a $k$-component dominating set of order less than $n/2-p$.
Hence, $D$ has the desired properties.
$\Box$ 

\medskip

\noindent If $G$, $u$, and $D$ are as in Lemma \ref{lemma5}(iv),
then $D$ is a {\it semi-$k$-component dominating set of $G$ with $u$ in the small component}.

\begin{lemma}[Shermer \cite{sh}]\label{lemma6}
Let $s$ and $n$ be positive integers with $s\geq 2$, and let $G$ be a maximal outerplanar graph of order $n$.
If $n\geq 2s$, then $G$ has a chord $xy$ such that one of the subgraphs of $G$ generated by $xy$
has $m$ $C(G)$-edges where $s\leq m\leq 2s-2$.
\end{lemma}
For a proof of Lemma \ref{lemma6}, the reader may refer to \cite{aldara,sh}.

\begin{lemma}\label{lemma7}
Let $k$ and $\ell$ be positive integers with $2\leq \ell\leq k$.
If $G$ is a maximal outerplanar graph of order $n=4k+2\ell$ that does not belong to ${\cal H}_k$, 
then $\gamma_k(G)\leq 2k+\ell-2=\left\lfloor\frac{kn}{2k+1}\right\rfloor$.
\end{lemma}
{\it Proof:} Note that $2k+\ell-2=\left\lfloor\frac{kn}{2k+1}\right\rfloor$
is equivalent to $0=\left\lfloor\frac{2k-\ell+2}{2k+1}\right\rfloor$,
which is equivalent to $0\leq 2k-\ell+2\leq 2k$.
Therefore, the equality $2k+\ell-2=\left\lfloor\frac{kn}{2k+1}\right\rfloor$ follows easily from $n=4k+2\ell$ and $2\leq \ell\leq k$.

It remains to show $\gamma_k(G)\leq 2k+\ell-2$.
For a contradiction, suppose that $G$ is a graph of order $4k+2\ell$ 
that does not belong to ${\cal H}_k$ and satisfies $\gamma_k(G)>2k+\ell-2$.
Since $n\geq 2(2k+2)$, Lemma \ref{lemma6} implies that $G$ has a chord $xy$ 
such that one of the subgraphs of $G$ generated by $xy$, say $G_{xy}$, 
has $m$ $C(G)$-edges where $2k+2\leq m\leq 4k+2$.
We assume that $xy$ is chosen such that $m$ is smallest possible subject to these conditions.
Let $G_z$ be the subgraph of $G$ generated by $xy$ that is distinct from $G_{xy}$.

If $m=2k+2$, then $G_{xy}$ has order $2k+3$, and $G_z$ has odd order $n-(2k+1)$. 
By Lemma \ref{lemma2}, Lemma \ref{lemma5}(i), and the choice of $G$,
the graph $G_{xy}$ has a $k$-component dominating set of order $k$, and 
$G_z$ has a $k$-component dominating set of order at most $\left\lfloor\frac{k(n-(2k+1))}{2k+1}\right\rfloor=k+\ell-2$
whose union is a $k$-component dominating set of $G$ of order at most $2k+\ell-2$,
which is a contradiction.
Hence, $m>2k+2$.

Let $z$ be the vertex of $G_{xy}$ such that $xyz$ is a triangle of $G$.
Let $G_x$ be the subgraph of $G$ generated by $yz$ that does not contain $x$,
and
let $G_y$ be the subgraph of $G$ generated by $xz$ that does not contain $y$.
Let $G_x$ and $G_y$ have orders $\ell_x+1$ and $\ell_y+1$, respectively.
Let $G_z$ have order $\ell_z+1$.
Note that $m=\ell_x+\ell_y$ and $n=\ell_x+\ell_y+\ell_z$.
The choice of $xy$ and $m>2k+2$ imply $\ell_x,\ell_y\geq 2$.

We consider different cases.

\medskip

\pagebreak

\noindent {\bf Case 1} {\it $\ell_x$ and $\ell_y$ are both odd.}

\medskip

\noindent By Lemma \ref{lemma4}(ii), 
$G_x$ has a $(\ell_x-1)/2$-component dominating set $D_x$ of order $(\ell_x-1)/2$ that intersects $yz$, and
$G_y$ has a $(\ell_y-1)/2$-component dominating set $D_y$ of order $(\ell_y-1)/2$ that intersects $xz$.
Note that $G[D_x\cup D_y]$ is a connected graph of order at least $(\ell_x-1)/2+(\ell_y-1)/2-1=m/2-2\geq k$.
Since $m$ is even, the order of $G_z$ is odd.

Suppose that $D_x$ and $D_y$ both contain $z$, which implies $|D_x\cup D_y|=m/2-2$.
By Lemma \ref{lemma4}(i), $G_z$ has a $\ell_z/2$-component dominating set $D_z$ of order $\ell_z/2$ that contains $x$,
and $D_x\cup D_y\cup D_z$ is a $k$-component dominating set of $G$ of order at most
$m/2-2+\ell_z/2=n/2-2=2k+\ell-2$,
which is a contradiction.
Hence, $D_x$ and $D_y$ do not both contain $z$.
By symmetry, we may assume that $D_y$ contains $x$.
Again, by Lemma \ref{lemma4}(i), $G_z$ has a $\ell_z/2$-component dominating set $D_z$ of order $\ell_z/2$ that contains $x$,
and $D_x\cup D_y\cup D_z$ is a $k$-component dominating set of $G$ of order at most
$m/2-1+\ell_z/2-1=2k+\ell-2$,
which is a contradiction.

\medskip

\noindent {\bf Case 2} {\it $\ell_x$ is odd and $\ell_y$ is even.}

\medskip

\noindent By Lemma \ref{lemma4}(ii), 
$G_x$ has a $(\ell_x-1)/2$-component dominating set $D_x$ of order $(\ell_x-1)/2$ that intersects $yz$.
Suppose that $D_x$ contains $z$.
By Lemma \ref{lemma4}(i), 
$G_y$ has a $\ell_y/2$-component dominating set $D_y$ of order $\ell_y/2$ that contains $z$.
Note that $G[D_x\cup D_y]$ is a connected graph of order $(\ell_x-1)/2+\ell_y/2-1=(m-3)/2\geq k$.
Since $m$ is odd, the order of $G_z$ is even.
By Lemma \ref{lemma4}(ii),
$G_z$ has a $(\ell_z-1)/2$-component dominating set $D_z$ of order $(\ell_z-1)/2$ that intersects $xy$,
and $D_x\cup D_y\cup D_z$ is a $k$-component dominating set of $G$ of order at most $(m-3)/2+(\ell_z-1)/2=n/2-2=2k+\ell-2$,
which is a contradiction.
Hence, $D_x$ contains $y$ but not $z$.
By Lemma \ref{lemma4}(i), 
$G_y$ has a $\ell_y/2$-component dominating set $D_y$ of order $\ell_y/2$ that contains $x$.
Since $D_z$ contains $x$ or $y$, 
the set $D_x\cup D_y\cup D_z$ is a $k$-component dominating set of $G$ of order at most $(m-1)/2+(\ell_z-1)/2-1=2k+\ell-2$,
which is a contradiction.

\medskip

\noindent {\bf Case 3} {\it $\ell_x$ and $\ell_y$ are both even.}

\medskip

\noindent By Lemma \ref{lemma4}(i), 
$G_x$ has a $\ell_x/2$-component dominating set $D_x$ of order $\ell_x/2$ that contains $z$, and
$G_y$ has a $\ell_y/2$-component dominating set $D_y$ of order $\ell_y/2$ that contains $z$.
Note that $G[D_x\cup D_y]$ is a connected graph of order $\ell_x/2+\ell_y/2-1=m/2-1>k$.
Since $m$ is even, the order of $G_z$ is odd.
By the choice of $xy$, we have $\ell_x\leq 2k$ and $\ell_y\leq 2k$,
which implies $\ell_y\geq 4$ and $\ell_x\geq 4$. 
If $m=4k+2$, then the choice of $xy$ implies $\ell_x=\ell_y=2k+1$, which is a contradiction.
Hence $m\leq 4k$, which implies that $\ell_z\geq 4$.

Let $N_{C_0(G_x)}(y)=\{ y',z\}$ and $N_{C_0(G_x)}(z)=\{ y,z'\}$.

Suppose that $y$ and $z$ both have degree at least $3$ in $G_x$.
By Lemma \ref{lemma4}(iii), 
$G_x$ has a $\ell_x/2$-component dominating set $D_x$ of order $\ell_x/2$ that contains $y$ and $z$.
By Lemma \ref{lemma4}(i),
$G_y$ has a $\ell_y/2$-component dominating set $D_y$ of order $\ell_y/2$ that contains $z$,
and
$G_z$ has a $\ell_z/2$-component dominating set $D_z$ of order $\ell_z/2$ that contains $y$.
Since $\ell_x/2+\ell_y/2+\ell_z/2-2=2k+\ell-2$,
the set $D_x\cup D_y\cup D_z$ is a $k$-component dominating set of $G$ of order at most $2k+\ell-2$,
which is a contradiction.
Hence, by symmetry, we may assume that $y$ has degree $2$ in $G_x$,
which implies that $z$ is adjacent to $y'$.

Suppose that $z$ has degree at least $4$ in $G_x$.
By Lemma \ref{lemma4}(iv), 
$G_x-y$ has a $(\ell_x/2-1)$-component dominating set $D_x'$ of order $\ell_x/2-1$ that contains $z$.
Choosing $D_y$ and $D_z$ as above, 
we obtain that $D_x\cup D_y\cup D_z$ is a $k$-component dominating set of $G$ of order at most $2k+\ell-2$,
which is a contradiction.
Hence, we may assume that $z$ has degree $3$ in $G_x$.

By symmetry, we obtain that 
$\{ d_{G_x}(y),d_{G_x}(z)\}=\{ d_{G_y}(x),d_{G_y}(z)\}=\{ 2,3\}$.
Furthermore, since our argument did not use the fact that $\ell_x\leq 2k$,
we also obtain, by symmetry, that $\{ d_{G_z}(x),d_{G_z}(y)\}=\{ 2,3\}$.

Suppose that $\ell_z\geq 2k+2$.
By Lemma \ref{lemma4}(i), 
$G_x$ has a $\ell_x/2$-component dominating set $D_x$ of order $\ell_x/2$ that contains $z$,
$G_y$ has a $\ell_y/2$-component dominating set $D_y$ of order $\ell_y/2$ that contains $z$, and
$G^-_z=G_z-\{ x,y\}$ has a $(\ell_z/2-1)$-component dominating set $D^-_z$ of order $\ell_z/2-1\geq k$.
Now, $D_x\cup D_y\cup D_z^-$ is a $k$-component dominating set of $G$ of order $2k+\ell-2$,
which is a contradiction.
Hence, $\ell_z\leq 2k$.

Suppose that $(G_x,yz)$ does not lie in ${\cal G}_{\ell_x}$.
Since the order of $G_x$ is $\ell_x+1$, 
the definition of ${\cal G}_{\ell_x}$ implies the existence of a 
$(\ell_x/2-2)$-component dominating set $D^-_x$ of $G^-_x=G_x-\{ y,z\}$ that intersects $y'z'$.
By Lemma \ref{lemma4}(i), 
$G_y$ has a $\ell_y/2$-component dominating set $D_y$ of order $\ell_y/2$ that contains $z$, and
$G_z$ has a $\ell_z/2$-component dominating set $D_z$ of order $\ell_z/2$ that contains $y$.
Now, $D^-_x\cup D_y\cup D_z$ is a $k$-component dominating set of $G$ of order at most $2k+\ell-2$,
which is a contradiction.

By symmetry, we obtain that 
$(G_x,yz)\in {\cal G}_{\ell_x}$,
$(G_y,xz)\in {\cal G}_{\ell_y}$, and
$(G_z,xy)\in {\cal G}_{\ell_z}$, 
which implies the contradiction $G\in {\cal H}_k^1\subseteq {\cal H}_k$. $\Box$

\medskip

\noindent We proceed to the proof of our main result.

\medskip

\noindent {\it Proof of Theorem \ref{theorem1}:}
Suppose for a contradiction, that $G$ is a counterexample of minimum order $n$.
Lemma \ref{lemma3} implies $n\geq 4k+4$.
Lemma \ref{lemma5} implies $G\not\in {\cal H}_k$ and $\gamma_k(G)>\left\lfloor\frac{kn}{2k+1}\right\rfloor$.

\begin{claim}\label{claim1}
$n\mod (2k+1)=2\ell$ for some $\ell\in [k-1]$.
\end{claim}
{\it Proof of Claim \ref{claim1}:} 
Suppose for a contradiction that $n\mod (2k+1)\not\in \{ 2\ell:\ell\in [k-1]\}.$

Clearly, $G$ contains no two adjacent vertices of degree $2$. 
If $G$ does not contain
either 
two vertices $u$ and $v$ of degree $2$ at distance $2$
or two adjacent vertices $u$ and $v$ such that $u$ has degree $2$ and $v$ has degree $3$,
then removing form $G$ all vertices of degree $2$ results in a maximal outerplanar graph
of minimum degree at least $3$, which is a contradiction.
Hence, let $u$ and $v$ have the stated properties.

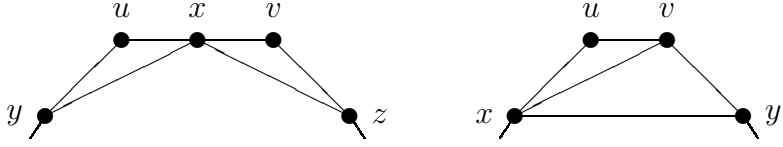
\begin{figure}[H]
\begin{center}
\unitlength 1mm 
\linethickness{0.4pt}
\ifx\plotpoint\undefined\newsavebox{\plotpoint}\fi 
\begin{picture}(49,19)(0,0)
\put(5,5){\circle*{2}}
\put(15,15){\circle*{2}}
\put(25,15){\circle*{2}}
\put(35,15){\circle*{2}}
\put(45,5){\circle*{2}}
\put(5,5){\line(1,1){10}}
\put(15,15){\line(1,0){20}}
\put(35,15){\line(1,-1){10}}
\put(45,5){\line(-2,1){20}}
\put(25,15){\line(-2,-1){20}}
\multiput(45,5)(.03333333,-.05){60}{\line(0,-1){.05}}
\multiput(5,5)(-.03333333,-.05){60}{\line(0,-1){.05}}
\put(1,5){\makebox(0,0)[cc]{$y$}}
\put(15,19){\makebox(0,0)[cc]{$u$}}
\put(25,19){\makebox(0,0)[cc]{$x$}}
\put(35,19){\makebox(0,0)[cc]{$v$}}
\put(49,5){\makebox(0,0)[cc]{$z$}}
\end{picture}\hspace{1cm}
\unitlength 1mm 
\linethickness{0.4pt}
\ifx\plotpoint\undefined\newsavebox{\plotpoint}\fi 
\begin{picture}(39,19)(0,0)
\put(5,5){\circle*{2}}
\put(15,15){\circle*{2}}
\put(25,15){\circle*{2}}
\put(5,5){\line(1,1){10}}
\put(25,15){\line(-2,-1){20}}
\multiput(5,5)(-.03333333,-.05){60}{\line(0,-1){.05}}
\put(1,5){\makebox(0,0)[cc]{$x$}}
\put(15,19){\makebox(0,0)[cc]{$u$}}
\put(25,19){\makebox(0,0)[cc]{$v$}}
\put(35,5){\circle*{2}}
\put(15,15){\line(1,0){10}}
\put(25,15){\line(1,-1){10}}
\put(35,5){\line(-1,0){30}}
\put(39,5){\makebox(0,0)[cc]{$y$}}
\multiput(35,5)(.03333333,-.05){60}{\line(0,-1){.05}}
\end{picture}
\end{center}
\caption{The two possibilities for $u$ and $v$ in the proof of the claim.}\label{fig4}
\end{figure}

The graph $G'=G-\{ u,v\}$ is a maximal outerplanar graph of order $n-2$.
In the first case, let $N_G(u)=\{ x,y\}$ and $N_G(v)=\{ x,z\}$.
Note that $xy$ and $yz$ are edges of $G'$ that belong to $C(G')$.
In the second case, let $N_G(u)=\{ v,x\}$ and $N_G(v)=\{ u,x,y\}$.
Note that $xy$ is an edge of $G'$ that belongs to $C(G')$.
See Figure \ref{fig4} for an illustration.

\noindent Suppose that $G'\not\in {\cal H}_k$.
By the choice of $G$, the graph $G'$ has a $k$-component dominating set $D'$ of order at most $\left\lfloor\frac{k(n-2)}{2k+1}\right\rfloor$,
and $D=D'\cup \{ x\}$ is a $k$-component dominating set of $G$ of order at most 
$\left\lfloor\frac{k(n-2)}{2k+1}\right\rfloor+1=\left\lfloor\frac{kn+1}{2k+1}\right\rfloor$.
Since $n\not\equiv 2\mod (2k+1)$, and $k$ and $2k+1$ are coprime,
we have $kn\not\equiv 2k\mod (2k+1)$, which implies the contradiction
$|D|\leq \left\lfloor\frac{kn+1}{2k+1}\right\rfloor=\left\lfloor\frac{kn}{2k+1}\right\rfloor$.
Hence, $G'\in {\cal H}_k^p$ for some $p\in [k-1]$.
By Lemma \ref{lemma5}(i), 
$n=n'+2$ is an even integer at least $4kp+2p+4=2p(2k+1)+4$ and at most $2k(2p+1)+2=2p(2k+1)+2(k-p+1)$,
which implies $n\mod (2k+1)\in \{ 2\ell:\ell\in [2,k-p+1]\}$.
Our assumption implies $n\mod (2k+1)=2k$, which implies $p=1$ and $n=2k(2p+1)+2=6k+2$.
If $G'$ arises as in the definition of ${\cal H}_k^1$ by suitably identifying vertices in three graphs 
from ${\cal G}_{\ell_1}$, ${\cal G}_{\ell_2}$, and ${\cal G}_{\ell_3}$, respectively,
then $\ell_i\leq 2k$ and $n'=6k$ imply $\ell_1=\ell_2=\ell_3=2k$.

If $x\not\in V(C_0(G'))$, then 
Lemma \ref{lemma5}(iv) implies that 
$G'$ has a semi-$k$-component dominating set $D'$ of order $\left\lfloor\frac{k(n-2)}{2k+1}\right\rfloor$ 
with $x$ in the small component, which is of order $k-1$.
In this case, let $D=D'\cup \{ u\}$.
If $x\in V(C_0(G'))$, then $y\not\in V(C_0(G'))$, and
Lemma \ref{lemma5}(iv) implies that 
$G'$ has a semi-$k$-component dominating set $D'$ of order $\left\lfloor\frac{k(n-2)}{2k+1}\right\rfloor$ 
with $x\not\in D'$ and $y$ in the small component, which is of order $k-1$.
In this case, let $D=D'\cup \{ x\}$.
In both cases, $D$ is a $k$-component dominating set of $G$ of order
$\left\lfloor\frac{kn+1}{2k+1}\right\rfloor=\left\lfloor\frac{kn}{2k+1}\right\rfloor$,
which is a contradiction. $\Box$

\medskip

\noindent Suppose $n\leq 6k+4$.
Since $n\geq 4k+4$, the claim implies $n=2(2k+1)+2\ell=4k+2(\ell+1)$ for some $\ell\in [k-1]$.
By Lemma \ref{lemma7}, $\gamma_k(G)\leq \left\lfloor\frac{kn}{2k+1}\right\rfloor$,
which is a contradiction.
Hence, $n\geq 6k+5$.

Since $6k+5\geq 2(2k+2)$, Lemma \ref{lemma4} implies the existence of a chord $xy$ 
such that one of the subgraphs of $G$ generated by $xy$, say $G_{xy}$, 
has $m$ $C(G)$-edges where $2k+2\leq m\leq 4k+2$.
We assume that $xy$ is chosen such that $m$ is smallest possible subject to these conditions.
Let $G_z$ denote the subgraph of $G$ generated by $xy$ distinct from $G_{xy}$.
Note that $G_z$ has order $n-m+1\geq (6k+5)-(4k+2)+1=2k+4$,
that is, contracting one or two edges of $G_z$ yields a graph of order at least $2k+2$.

Suppose $m=2k+2$.
If $G_z\not\in {\cal H}_k$, 
then Lemma \ref{lemma2} and the choice of $G$ imply
$$\gamma_k(G)\leq \gamma_k(G_{xy})+\gamma_k(G_z)\leq 
k+\left\lfloor\frac{k(n-(2k+1))}{2k+1}\right\rfloor
=\left\lfloor\frac{kn}{2k+1}\right\rfloor,$$
which is a contradiction.
Hence, $G_z\in {\cal H}_k$.
Since $xy$ is an edge of $C_0(G_z)$,
we obtain that either $x$ or $y$ does not belong to $C_0(G_z)$.
By symmetry, we may assume that $x$ does not belong to $C_0(G_z)$.
By Lemma \ref{lemma5}(iv),
$G_z$ has a semi-$k$-component dominating set $D_z$ of order $\left\lfloor\frac{k(n-(2k+1))}{2k+1}\right\rfloor$ 
with $x$ in the small component.
By Lemma \ref{lemma4}(i), $G_{xy}$ has a $(k+1)$-component dominating set $D_{xy}$ of order $k+1$ that contains $x$.
Now, $D_{xy}\cup D_z$ is a $k$-component dominating set of $G$ of order at most
$(k+1)+\left\lfloor\frac{k(n-(2k+1))}{2k+1}\right\rfloor-1=\left\lfloor\frac{kn}{2k+1}\right\rfloor$,
which is a contradiction.
Hence, $m>2k+2$.

Let $z$ be the vertex of $G_{xy}$ such that $xyz$ is a triangle of $G$.
Let $G_x$ be the subgraph of $G$ generated by $yz$ that does not contain $x$,
and
let $G_y$ be the subgraph of $G$ generated by $xz$ that does not contain $y$.
Let $G_x$ and $G_y$ have orders $\ell_x+1$ and $\ell_y+1$, respectively.
The choice of $xy$ and $m>2k+2$ imply $\ell_x,\ell_y\geq 2$.
Let $G_z$ have order $\ell_z+1$.
Note that $m=\ell_x+\ell_y$ and $n=\ell_x+\ell_y+\ell_z$.

We consider different cases.

\medskip

\noindent {\bf Case 1} {\it $\ell_x$ and $\ell_y$ are both odd.}

\medskip

\noindent By Lemma \ref{lemma4}(i), 
$G_x$ has a $(\ell_x-1)/2$-component dominating set $D_x$ of order $(\ell_x-1)/2$ that intersects $yz$, and
$G_y$ has a $(\ell_y-1)/2$-component dominating set $D_y$ of order $(\ell_y-1)/2$ that intersects $xz$.
Since $m\geq 2k+4$, $G[D_x\cup D_y]$ is a connected graph of order at least $m/2-2\geq k$.

\medskip

\noindent {\bf Case 1.1} {\it $z\in D_x\cap D_y$.}

\medskip

\noindent Suppose $G_z\not\in {\cal H}_k$.
We obtain that $\gamma_k(G)$ is at most
\begin{eqnarray*}
|D_x\cup D_y|+\gamma_k(G_z)
\leq m/2-2+\left\lfloor\frac{k(n-m+1)}{2k+1}\right\rfloor
= \left\lfloor\frac{kn}{2k+1}+\frac{m/2-3k-2}{2k+1}\right\rfloor
\leq \left\lfloor\frac{kn}{2k+1}\right\rfloor,
\end{eqnarray*}
where the last inequality is implied by $m<6k+4$,
which is a contradiction.
Hence, $G_z\in {\cal H}_k$.
By symmetry, we may assume that $x$ does not belong to $C_0(G_z)$.
By Lemma \ref{lemma5}(iv),
$G_z$ has a semi-$k$-component dominating set $D_z$ of order $\left\lfloor\frac{k(n-m+1)}{2k+1}\right\rfloor$ 
with $x$ in the small component.
Again,
$\gamma_k(G)\leq m/2-2+\left\lfloor\frac{k(n-m+1)}{2k+1}\right\rfloor
\leq \left\lfloor\frac{kn}{2k+1}\right\rfloor$,
which is a contradiction.

\medskip

\noindent {\bf Case 1.2} {\it $z\not\in D_x\cap D_y$.}

\medskip

\noindent By symmetry, we may assume that $y\in D_x$.
Let $G_z'$ arise from $G_z$ by contracting the edge $xy$ to a new vertex $u^*$.
Suppose $G'_z\not\in {\cal H}_k$.
By the choice of $G$, the graph $G'_z$ has a $k$-component dominating set $D_z'$ of order at most $\left\lfloor\frac{k(n-m)}{2k+1}\right\rfloor$.
If $u^*\not\in D_z'$, then let $D=D_x\cup D_y\cup D_z'$.
If $u^*\in D_z'$, then let $D=D_x\cup D_y\cup (D_z'\setminus \{ u^*\})\cup \{ x\}$.
In both cases, $D$ is a $k$-component dominating set of $G$ of order at most 
\begin{eqnarray*}
m/2-1+\left\lfloor\frac{k(n-m)}{2k+1}\right\rfloor
& = & \left\lfloor\frac{kn}{2k+1}+\frac{m/2-2k-1}{2k+1}\right\rfloor
\leq \left\lfloor\frac{kn}{2k+1}\right\rfloor,
\end{eqnarray*}
where the last inequality is implied by $m\leq 4k+2$,
which is a contradiction.
Hence, $G'_z\in {\cal H}_k$.
Suppose $u^*\not\in C_0(G_z')$.
By Lemma \ref{lemma5}(iv), 
$G'_z$ has a semi-$k$-component dominating set $D'_z$ of order $\left\lfloor\frac{k(n-m)}{2k+1}\right\rfloor$ 
with $u^*$ in the small component.
The set $D_x\cup D_y\cup (D_z'\setminus \{ u^*\})\cup \{ x\}$
is a $k$-component dominating set of $G$ of order at most 
$m/2-1+\left\lfloor\frac{k(n-m)}{2k+1}\right\rfloor\leq \left\lfloor\frac{kn}{2k+1}\right\rfloor$,
which is a contradiction.
Hence, $u^*\in C_0(G'_z)$.
If $N_{C(G_z)}(y)=\{ x,y'\}$, then $y'\not\in C_0(G'_z)$.
By Lemma \ref{lemma5}(iv), 
$G'_z$ has a semi-$k$-component dominating set $D'_z$ of order $\left\lfloor\frac{k(n-m)}{2k+1}\right\rfloor$ 
with $y'$ in the small component such that $u^*\not\in D_z'$.
The set $D_x\cup D_y\cup D_z'$
is a $k$-component dominating set of $G$ of order at most 
$m/2-1+\left\lfloor\frac{k(n-m)}{2k+1}\right\rfloor\leq \left\lfloor\frac{kn}{2k+1}\right\rfloor$,
which is a contradiction.

\medskip

\pagebreak

\noindent {\bf Case 2} {\it $\ell_x$ is odd and $\ell_y$ is even.}

\medskip

\noindent Note that $m$ is odd, which implies $m\leq 4k+1$.
By Lemma \ref{lemma4}(ii), 
$G_x$ has a $(\ell_x-1)/2$-component dominating set $D_x$ of order $(\ell_x-1)/2$ that intersects $yz$.

\medskip

\noindent {\bf Case 2.1} {\it $z\in D_x$.}

\medskip

\noindent By Lemma \ref{lemma4}(i), 
$G_y$ has a $\ell_y/2$-component dominating set $D_y$ of order $\ell_y/2$ that contains $z$.
Since $m\geq 2k+3$, we have $|D_x\cup D_y|=m/2-3/2\geq k$.
Suppose $G_z\not\in {\cal H}_k$.
By the choice of $G$, the graph $G_z$ has a $k$-component dominating set $D_z$ of order at most $\left\lfloor\frac{k(n-m+1)}{2k+1}\right\rfloor$.
Now, the set 
$D_x\cup D_y\cup D_z$ is a $k$-component dominating set of $G$ of order at most 
\begin{eqnarray*}
m/2-3/2+\left\lfloor\frac{k(n-m+1)}{2k+1}\right\rfloor
& = & \left\lfloor\frac{kn}{2k+1}+\frac{m/2-2k-3/2}{2k+1}\right\rfloor
\leq \left\lfloor\frac{kn}{2k+1}\right\rfloor,
\end{eqnarray*}
where the last inequality is implied by $m<4k+3$,
which is a contradiction.
Hence, $G_z\in {\cal H}_k$.
By Lemma \ref{lemma5}(iv), 
$G_z$ has a semi-$k$-component dominating set $D_z$ of order $\left\lfloor\frac{k(n-m+1)}{2k+1}\right\rfloor$ 
with $x$ or $y$ in the small component.
Again, the set
$D_x\cup D_y\cup D_z$ is a $k$-component dominating set of $G$ of order at most 
$m/2-3/2+\left\lfloor\frac{k(n-m+1)}{2k+1}\right\rfloor\leq \left\lfloor\frac{kn}{2k+1}\right\rfloor$,
which is a contradiction.

\medskip

\noindent {\bf Case 2.2} {\it $z\not\in D_x$.}

\medskip

\noindent We have $y\in D_x$.
By Lemma \ref{lemma4}(i), 
$G_y$ has a $\ell_y/2$-component dominating set $D_y$ of order $\ell_y/2$ that contains $x$.
Note that $|D_x\cup D_y|=m/2-1/2>k$.
Let $N_{C(G_z)}(y)=\{ x,y'\}$ and $N_{C(G_z)}(x)=\{ x',y\}$.
Let $G_z''$ arise from $G_z$ by contracting the two edges $xy$ and $yy'$ to a new vertex $u^*$.
Suppose $G''_z\not\in {\cal H}_k$.
By the choice of $G$, the graph $G''_z$ has a $k$-component dominating set $D''_z$ 
of order at most $\left\lfloor\frac{k(n-m-1)}{2k+1}\right\rfloor$.
If $u^*\in D_z''$, then let $D=D_x\cup D_y\cup (D_z''\setminus \{ u^*\})\cup \{ y'\}$.
If $u^*\not\in D_z''$, then let $D=D_x\cup D_y\cup D_z''$.
The set $D$ is a $k$-component dominating set o $G$ of order at most 
\begin{eqnarray*}
m/2-1/2+\left\lfloor\frac{k(n-m-1)}{2k+1}\right\rfloor
& = & \left\lfloor\frac{kn}{2k+1}+\frac{m/2-2k-1/2}{2k+1}\right\rfloor
\leq \left\lfloor\frac{kn}{2k+1}\right\rfloor
\end{eqnarray*}
where the last inequality is implied by $m\leq 4k+1$,
which is a contradiction.
Hence, $G''_z\in {\cal H}_k$.
Suppose $u^*\not\in C_0(G_z'')$.
By Lemma \ref{lemma5}(iv), 
$G''_z$ has a semi-$k$-component dominating set $D''_z$ of order $\left\lfloor\frac{k(n-m-1)}{2k+1}\right\rfloor$ 
with $u^*$ in the small component.
The set $D_x\cup D_y\cup (D_z''\setminus \{ u^*\})\cup \{ y'\}$
is a $k$-component dominating set of $G$ of order at most 
$m/2-1/2+\left\lfloor\frac{k(n-m-1)}{2k+1}\right\rfloor\leq \left\lfloor\frac{kn}{2k+1}\right\rfloor$,
which is a contradiction.
Hence, $u^*\in C_0(G'_z)$.
By Lemma \ref{lemma5}(iv), 
$G''_z$ has a semi-$k$-component dominating set $D''_z$ of order $\left\lfloor\frac{k(n-m-1)}{2k+1}\right\rfloor$ 
with $x'$ in the small component such that $u^*\not\in D_z''$.
The set $D_x\cup D_y\cup D_z''$
is a $k$-component dominating set of $G$ of order at most 
$m/2-1/2+\left\lfloor\frac{k(n-m-1)}{2k+1}\right\rfloor\leq \left\lfloor\frac{kn}{2k+1}\right\rfloor$,
which is a contradiction.

\medskip

\noindent {\bf Case 3} {\it $\ell_x$ and $\ell_y$ are both even.}

\medskip

\noindent We have $m\geq 2k+4$.
By the choice of $xy$, this implies $4\leq \ell_x,\ell_y\leq 2k$, and, hence, $m\leq 4k$.
Let $G_z'$ arise from $G_z$ by contracting the edge $xy$ to a new vertex $u^*$.

Suppose $G'_z\not\in {\cal H}_k$.
By the choice of $G$, 
the graph $G'_z$ has a $k$-component dominating set $D'_z$ of order at most $\left\lfloor\frac{k(n-m)}{2k+1}\right\rfloor$.
If $u^*\in D_z'$, then Lemma \ref{lemma4}(i) implies that 
$G_x$ has a $\ell_x/2$-component dominating set $D_x$ of order $\ell_x/2$ that contains $y$, and 
$G_y$ has a $\ell_y/2$-component dominating set $D_y$ of order $\ell_y/2$ that contains $x$.
In this case, let $D=D_x\cup D_y\cup (D_z'\setminus \{ u^*\})$.
If $u^*\not\in D_z'$, then Lemma \ref{lemma4}(i) implies that 
$G_x$ has a $\ell_x/2$-component dominating set $D_x$ of order $\ell_x/2$ that contains $z$, and 
$G_y$ has a $\ell_y/2$-component dominating set $D_y$ of order $\ell_y/2$ that contains $z$.
In this case, let $D=D_x\cup D_y\cup D_z'$.
In both cases $D$ is a $k$-component dominating set of $G$ of order at most
$m/2-1+\left\lfloor\frac{k(n-m)}{2k+1}\right\rfloor\leq \left\lfloor\frac{kn}{2k+1}\right\rfloor$
(cf. the calculation in Case 1.2), which is a contradiction.
Hence, $G'_z\in {\cal H}_k$.

Suppose $u^*\not\in C_0(G_z')$.
By Lemma \ref{lemma4}(i),
$G_x$ has a $\ell_x/2$-component dominating set $D_x$ of order $\ell_x/2$ that contains $y$, and 
$G_y$ has a $\ell_y/2$-component dominating set $D_y$ of order $\ell_y/2$ that contains $x$.
Note that $k<m/2-1\leq |D_x\cup D_y|\leq m/2$.
By Lemma \ref{lemma5}(iv), 
$G'_z$ has a semi-$k$-component dominating set $D'_z$ of order $\left\lfloor\frac{k(n-m)}{2k+1}\right\rfloor$ 
with $u^*$ in the small component.
The set $D_x\cup D_y\cup (D'_z\setminus \{ u^*\})$
is a $k$-component dominating set of $G$ of order at most 
$m/2-1+\left\lfloor\frac{k(n-m)}{2k+1}\right\rfloor\leq \left\lfloor\frac{kn}{2k+1}\right\rfloor$,
which is a contradiction.
Hence, $u^*\in C_0(G_z')$.

Let $C_0(G_z')$ have order $2p+1$ for some $p\in [k-1]$.
For $i\in [2p+1]$, let $\ell_i$ with $4\leq \ell_i\leq 2k$
and $(G_i,x_iy_i)\in {\cal G}_{\ell_i}$ 
be as in the definition of ${\cal H}_k^p$
such that $G_z'$ arises by suitably identifying vertices in the graphs $G_1,\ldots,G_{2p+1}$,
that is, $C_0(G_z')$ is the cycle $x_1x_2\ldots x_{2p+1}x_1$.
Let $u^*=x_1$, that is, $xzyx_2\ldots x_{2p+1}x$ is a cycle in $G$.
For $i\in [2p+1]$, let $D_i^x$ and $D_i^y$ be as in the proof of Lemma \ref{lemma5}(iii).

By Lemma \ref{lemma4}(i),
$G_x$ has a $\ell_x/2$-component dominating set $D_x$ of order $\ell_x/2$ that contains $z$, and 
$G_y$ has a $\ell_y/2$-component dominating set $D_y$ of order $\ell_y/2$ that contains $z$.
Now, the set
$$D_x\cup D_y\cup D_1^y\cup D_2^x\cup D_3^y\cup D_4^x\cup\cdots\cup D_{2p-1}^y\cup D_{2p}^x\cup D_{2p+1}^y$$
is a $k$-component dominating set of $G$ of order at most $n/2-(p+1)$.
By the choice of $G$, we obtain $\lfloor\frac{kn}{2k+1}\rfloor<\gamma_k(G)\leq n/2-(p+1)$,
which implies $n\geq 4k(p+1)+2(p+1)+2$.

Suppose that $p=k-1$.
We obtain $n\geq 4k(p+1)+2(p+1)+2=4k^2+2k+2$
as well as $n\leq n(G_z')+m\leq (2k-1)2k+m\leq 4k^2+2k$,
which is a contradiction.
Hence, $p\leq k-2$, which implies $k\geq 3$.

Suppose that $d_{G_y}(x),d_{G_y}(z)\geq 3$.
By Lemma \ref{lemma4}(iii),
$G_y$ has a $\ell_y/2$-component dominating set $D_y$ of order $\ell_y/2$ that contains $x$ and $z$.
Choosing $D_x$ as above, the set
$$D_x\cup D_y\cup D_1^y\cup D_2^x\cup D_3^y\cup D_4^x\cup\cdots\cup D_{2p-1}^y\cup D_{2p}^x\cup D_{2p+1}^y$$
is a $k$-component dominating set of $G$ of order at most $n/2-(p+2)$.
Since $n\leq (2p+1)2k+2k$, we obtain $\gamma_k(G)\leq n/2-(p+2)\leq \left\lfloor\frac{kn}{2k+1}\right\rfloor$,
which is a contradiction. 
Hence, one of the two degrees $d_{G_y}(x)$ and $d_{G_y}(z)$ is $2$.

Suppose that $d_{G_y}(x)=2$ and $d_{G_y}(z)\geq 4$.
By Lemma \ref{lemma4}(iv),
$G_y-x$ has a $(\ell_y/2-1)$-component dominating set $D_y$ of order $\ell_y/2-1$ that contains $z$.
Choosing $D_x$ as above, the set
$$D_x\cup D_y\cup D_1^y\cup D_2^x\cup D_3^y\cup D_4^x\cup\cdots\cup D_{2p-1}^y\cup D_{2p}^x\cup D_{2p+1}^y$$
is a $k$-component dominating set of $G$ of order at most $n/2-(p+2)\leq \left\lfloor\frac{kn}{2k+1}\right\rfloor$,
which is a contradiction. 
Hence, if $d_{G_y}(x)=2$, then $d_{G_y}(z)=3$.
By a symmetric argument, we obtain that if $d_{G_y}(z)=2$, then $d_{G_y}(x)=3$,
that is, $\{ d_{G_y}(x),d_{G_y}(z)\}=\{ 2,3\}$.
By symmetry, $\{ d_{G_x}(y),d_{G_x}(z)\}=\{ 2,3\}$.

Let $N_{C(G_y)}(x)=\{ x',z\}$ and $N_{C(G_y)}(z)=\{ x,z'\}$.

Suppose that $(G_y,xz)$ does not belong to ${\cal G}_{\ell_y}$.
By the definition of ${\cal G}_{\ell_y}$, 
this implies that the graph $G_y^-=G_y-\{ x,z\}$ has a $(\ell_y/2-2)$-component dominating set $D_y^-$ of order $\ell_y/2-2$ that intersects $x'z'$.
Let $D_x$ be as above.
Now, the set 
$$D_x\cup D^-_y\cup D_1^y\cup D_2^x\cup D_3^y\cup D_4^x\cup\cdots\cup D_{2p-1}^y\cup D_{2p}^x\cup D_{2p+1}^y$$
contains $x\in D_{2p+1}^y$ and $z\in D_x$,
which implies that it is a $k$-component dominating set of $G$ of order at most $n/2-(p+2)\leq \left\lfloor\frac{kn}{2k+1}\right\rfloor$,
which is a contradiction. 
Hence, $(G_y,xz)\in {\cal G}_{\ell_y}$, and, by symmetry, $(G_x,yz)\in {\cal G}_{\ell_x}$.
Altogether, this implies that $G\in {\cal H}_k^{p+1}$,
which is the final contradiction and completes the proof. $\Box$

\medskip

\begin{figure}[H]
\begin{center}
\unitlength 0.6mm 
\linethickness{0.4pt}
\ifx\plotpoint\undefined\newsavebox{\plotpoint}\fi 
\begin{picture}(12,52)(0,0)
\put(10,10){\circle*{2}}
\put(10,20){\circle*{2}}
\put(10,30){\circle*{2}}
\put(10,40){\circle*{2}}
\put(10,50){\circle*{2}}
\put(0,10){\circle*{2}}
\put(0,20){\circle*{2}}
\put(0,30){\circle*{2}}
\put(0,40){\circle*{2}}
\put(0,50){\circle*{2}}
\put(10,10){\line(-1,0){10}}
\put(10,20){\line(-1,0){10}}
\put(10,30){\line(-1,0){10}}
\put(10,40){\line(-1,0){10}}
\put(10,50){\line(-1,0){10}}
\put(0,10){\line(0,1){10}}
\put(0,20){\line(0,1){10}}
\put(0,30){\line(0,1){10}}
\put(0,40){\line(0,1){10}}
\put(0,20){\line(1,-1){10}}
\put(0,30){\line(1,-1){10}}
\put(0,40){\line(1,-1){10}}
\put(0,50){\line(1,-1){10}}
\put(10,10){\line(0,1){8}}
\put(10,20){\line(0,1){8}}
\put(10,30){\line(0,1){8}}
\put(10,40){\line(0,1){8}}
\put(10,18){\line(0,1){2}}
\put(10,28){\line(0,1){2}}
\put(10,38){\line(0,1){2}}
\put(10,48){\line(0,1){2}}
\put(10,10){\circle{4}}
\put(10,20){\circle{4}}
\put(10,30){\circle{4}}
\put(10,40){\circle{4}}
\put(10,50){\circle{4}}
\put(5,0){\circle*{2}}
\put(10,10){\line(-1,-2){5}}
\put(5,0){\line(-1,2){5}}
\end{picture}\hspace{1cm}
\linethickness{0.4pt}
\ifx\plotpoint\undefined\newsavebox{\plotpoint}\fi 
\begin{picture}(12,111)(0,0)
\put(10,10){\circle*{2}}
\put(10,20){\circle*{2}}
\put(10,30){\circle*{2}}
\put(10,40){\circle*{2}}
\put(10,50){\circle*{2}}
\put(10,60){\circle*{2}}
\put(10,70){\circle*{2}}
\put(10,80){\circle*{2}}
\put(10,90){\circle*{2}}
\put(10,100){\circle*{2}}
\put(0,10){\circle*{2}}
\put(0,20){\circle*{2}}
\put(0,30){\circle*{2}}
\put(0,40){\circle*{2}}
\put(0,50){\circle*{2}}
\put(0,60){\circle*{2}}
\put(0,70){\circle*{2}}
\put(0,80){\circle*{2}}
\put(0,90){\circle*{2}}
\put(0,100){\circle*{2}}
\put(10,10){\line(-1,0){10}}
\put(10,20){\line(-1,0){10}}
\put(10,30){\line(-1,0){10}}
\put(10,40){\line(-1,0){10}}
\put(10,50){\line(-1,0){10}}
\put(10,60){\line(-1,0){10}}
\put(10,70){\line(-1,0){10}}
\put(10,80){\line(-1,0){10}}
\put(10,90){\line(-1,0){10}}
\put(10,100){\line(-1,0){10}}
\put(0,10){\line(0,1){10}}
\put(0,20){\line(0,1){10}}
\put(0,30){\line(0,1){10}}
\put(0,40){\line(0,1){10}}
\put(0,50){\line(0,1){10}}
\put(0,60){\line(0,1){10}}
\put(0,70){\line(0,1){10}}
\put(0,80){\line(0,1){10}}
\put(0,90){\line(0,1){10}}
\put(0,20){\line(1,-1){10}}
\put(0,30){\line(1,-1){10}}
\put(0,40){\line(1,-1){10}}
\put(0,50){\line(1,-1){10}}
\put(0,60){\line(1,-1){10}}
\put(0,70){\line(1,-1){10}}
\put(0,80){\line(1,-1){10}}
\put(0,90){\line(1,-1){10}}
\put(0,100){\line(1,-1){10}}
\put(10,10){\line(0,1){8}}
\put(10,20){\line(0,1){8}}
\put(10,30){\line(0,1){8}}
\put(10,40){\line(0,1){8}}
\put(10,50){\line(0,1){8}}
\put(10,60){\line(0,1){8}}
\put(10,70){\line(0,1){8}}
\put(10,80){\line(0,1){8}}
\put(10,90){\line(0,1){8}}
\put(10,18){\line(0,1){2}}
\put(10,28){\line(0,1){2}}
\put(10,38){\line(0,1){2}}
\put(10,48){\line(0,1){2}}
\put(10,58){\line(0,1){2}}
\put(10,68){\line(0,1){2}}
\put(10,78){\line(0,1){2}}
\put(10,88){\line(0,1){2}}
\put(10,98){\line(0,1){2}}
\put(10,10){\circle{4}}
\put(10,20){\circle{4}}
\put(10,30){\circle{4}}
\put(10,40){\circle{4}}
\put(10,50){\circle{4}}
\put(10,60){\circle{4}}
\put(10,70){\circle{4}}
\put(10,80){\circle{4}}
\put(10,90){\circle{4}}
\put(10,100){\circle{4}}
\put(-1,60){\line(1,0){1}}
\put(5,0){\circle*{2}}
\put(0,10){\line(1,-2){5}}
\put(5,0){\line(1,2){5}}
\put(5,110){\circle*{2}}
\put(10,100){\line(-1,2){5}}
\put(5,110){\line(-1,-2){5}}
\end{picture}\hspace{1cm}
\linethickness{0.4pt}
\ifx\plotpoint\undefined\newsavebox{\plotpoint}\fi 
\begin{picture}(111,91)(0,0)
\put(10,0){\circle*{2}}
\put(10,20){\circle*{2}}
\put(20,0){\circle*{2}}
\put(20,20){\circle*{2}}
\put(30,0){\circle*{2}}
\put(30,20){\circle*{2}}
\put(40,0){\circle*{2}}
\put(40,20){\circle*{2}}
\put(50,0){\circle*{2}}
\put(50,20){\circle*{2}}
\put(60,0){\circle*{2}}
\put(60,20){\circle*{2}}
\put(70,0){\circle*{2}}
\put(70,20){\circle*{2}}
\put(80,0){\circle*{2}}
\put(80,20){\circle*{2}}
\put(60,50){\circle*{2}}
\put(90,0){\circle*{2}}
\put(90,20){\circle*{2}}
\put(60,60){\circle*{2}}
\put(60,70){\circle*{2}}
\put(100,0){\circle*{2}}
\put(100,20){\circle*{2}}
\put(60,80){\circle*{2}}
\put(10,10){\circle*{2}}
\put(10,30){\circle*{2}}
\put(20,10){\circle*{2}}
\put(20,30){\circle*{2}}
\put(30,10){\circle*{2}}
\put(30,30){\circle*{2}}
\put(40,10){\circle*{2}}
\put(40,30){\circle*{2}}
\put(50,10){\circle*{2}}
\put(50,30){\circle*{2}}
\put(60,10){\circle*{2}}
\put(70,10){\circle*{2}}
\put(70,30){\circle*{2}}
\put(80,10){\circle*{2}}
\put(80,30){\circle*{2}}
\put(50,50){\circle*{2}}
\put(90,10){\circle*{2}}
\put(90,30){\circle*{2}}
\put(50,60){\circle*{2}}
\put(50,70){\circle*{2}}
\put(100,10){\circle*{2}}
\put(100,30){\circle*{2}}
\put(50,80){\circle*{2}}
\put(10,0){\line(0,1){10}}
\put(10,20){\line(0,1){10}}
\put(20,0){\line(0,1){10}}
\put(20,20){\line(0,1){10}}
\put(30,0){\line(0,1){10}}
\put(30,20){\line(0,1){10}}
\put(40,0){\line(0,1){10}}
\put(40,20){\line(0,1){10}}
\put(50,0){\line(0,1){10}}
\put(50,20){\line(0,1){10}}
\put(60,0){\line(0,1){10}}
\put(60,20){\line(0,1){10}}
\put(70,0){\line(0,1){10}}
\put(70,20){\line(0,1){10}}
\put(80,0){\line(0,1){10}}
\put(80,20){\line(0,1){10}}
\put(60,50){\line(-1,0){10}}
\put(90,0){\line(0,1){10}}
\put(90,20){\line(0,1){10}}
\put(60,60){\line(-1,0){10}}
\put(60,70){\line(-1,0){10}}
\put(100,0){\line(0,1){10}}
\put(100,20){\line(0,1){10}}
\put(60,80){\line(-1,0){10}}
\put(10,10){\line(1,0){10}}
\put(10,30){\line(1,0){10}}
\put(20,10){\line(1,0){10}}
\put(20,30){\line(1,0){10}}
\put(30,10){\line(1,0){10}}
\put(30,30){\line(1,0){10}}
\put(40,10){\line(1,0){10}}
\put(40,30){\line(1,0){10}}
\put(50,10){\line(1,0){10}}
\put(50,30){\line(1,0){10}}
\put(60,10){\line(1,0){10}}
\put(60,30){\line(1,0){10}}
\put(70,10){\line(1,0){10}}
\put(70,30){\line(1,0){10}}
\put(50,40){\line(0,1){10}}
\put(80,10){\line(1,0){10}}
\put(80,30){\line(1,0){10}}
\put(50,50){\line(0,1){10}}
\put(90,10){\line(1,0){10}}
\put(90,30){\line(1,0){10}}
\put(50,60){\line(0,1){10}}
\put(50,70){\line(0,1){10}}
\put(20,10){\line(-1,-1){10}}
\put(20,30){\line(-1,-1){10}}
\put(30,10){\line(-1,-1){10}}
\put(30,30){\line(-1,-1){10}}
\put(40,10){\line(-1,-1){10}}
\put(40,30){\line(-1,-1){10}}
\put(50,10){\line(-1,-1){10}}
\put(50,30){\line(-1,-1){10}}
\put(60,10){\line(-1,-1){10}}
\put(60,30){\line(-1,-1){10}}
\put(70,10){\line(-1,-1){10}}
\put(70,30){\line(-1,-1){10}}
\put(80,10){\line(-1,-1){10}}
\put(80,30){\line(-1,-1){10}}
\put(50,50){\line(1,-1){10}}
\put(90,10){\line(-1,-1){10}}
\put(90,30){\line(-1,-1){10}}
\put(50,60){\line(1,-1){10}}
\put(100,10){\line(-1,-1){10}}
\put(100,30){\line(-1,-1){10}}
\put(50,70){\line(1,-1){10}}
\put(50,80){\line(1,-1){10}}
\put(10,0){\line(1,0){8}}
\put(10,20){\line(1,0){8}}
\put(20,0){\line(1,0){8}}
\put(20,20){\line(1,0){8}}
\put(30,0){\line(1,0){8}}
\put(30,20){\line(1,0){8}}
\put(40,0){\line(1,0){8}}
\put(40,20){\line(1,0){8}}
\put(50,0){\line(1,0){8}}
\put(50,20){\line(1,0){8}}
\put(60,0){\line(1,0){8}}
\put(60,20){\line(1,0){8}}
\put(70,0){\line(1,0){8}}
\put(70,20){\line(1,0){8}}
\put(60,40){\line(0,1){8}}
\put(80,0){\line(1,0){8}}
\put(80,20){\line(1,0){8}}
\put(60,50){\line(0,1){8}}
\put(90,0){\line(1,0){8}}
\put(90,20){\line(1,0){8}}
\put(60,60){\line(0,1){8}}
\put(60,70){\line(0,1){8}}
\put(18,0){\line(1,0){2}}
\put(18,20){\line(1,0){2}}
\put(28,0){\line(1,0){2}}
\put(28,20){\line(1,0){2}}
\put(38,0){\line(1,0){2}}
\put(38,20){\line(1,0){2}}
\put(48,0){\line(1,0){2}}
\put(48,20){\line(1,0){2}}
\put(58,0){\line(1,0){2}}
\put(58,20){\line(1,0){2}}
\put(60,38){\line(0,1){2}}
\put(68,0){\line(1,0){2}}
\put(68,20){\line(1,0){2}}
\put(78,0){\line(1,0){2}}
\put(78,20){\line(1,0){2}}
\put(60,48){\line(0,1){2}}
\put(88,0){\line(1,0){2}}
\put(88,20){\line(1,0){2}}
\put(60,58){\line(0,1){2}}
\put(60,68){\line(0,1){2}}
\put(98,0){\line(1,0){2}}
\put(98,20){\line(1,0){2}}
\put(60,78){\line(0,1){2}}
\put(10,0){\circle{4}}
\put(10,20){\circle{4}}
\put(20,0){\circle{4}}
\put(20,20){\circle{4}}
\put(30,0){\circle{4}}
\put(30,20){\circle{4}}
\put(40,0){\circle{4}}
\put(40,20){\circle{4}}
\put(50,0){\circle{4}}
\put(50,20){\circle{4}}
\put(60,0){\circle{4}}
\put(60,20){\circle{4}}
\put(60,40){\circle{4}}
\put(70,0){\circle{4}}
\put(70,20){\circle{4}}
\put(80,0){\circle{4}}
\put(80,20){\circle{4}}
\put(60,50){\circle{4}}
\put(90,0){\circle{4}}
\put(90,20){\circle{4}}
\put(60,60){\circle{4}}
\put(60,70){\circle{4}}
\put(100,0){\circle{4}}
\put(100,20){\circle{4}}
\put(60,80){\circle{4}}
\put(50,20){\line(0,-1){9}}
\put(60,20){\line(0,-1){9}}
\put(60,11){\line(0,-1){1}}
\put(60,20){\line(-1,-1){10}}
\put(60,40){\line(0,-1){11}}
\put(60,30){\line(-1,1){10}}
\put(50,40){\line(0,-1){10}}
\put(60,40){\circle*{2}}
\put(50,40){\circle*{2}}
\put(50,40){\line(1,0){10}}
\put(60,30){\circle*{2}}
\put(0,25){\circle*{2}}
\put(0,5){\circle*{2}}
\put(110,5){\circle*{2}}
\put(110,25){\circle*{2}}
\put(100,30){\line(2,-1){10}}
\put(110,25){\line(-2,-1){10}}
\put(100,10){\line(2,-1){10}}
\put(110,5){\line(-2,-1){10}}
\put(10,30){\line(-2,-1){10}}
\put(0,25){\line(2,-1){10}}
\put(10,10){\line(-2,-1){10}}
\put(0,5){\line(2,-1){10}}
\put(55,90){\circle*{2}}
\put(60,80){\line(-1,2){5}}
\put(55,90){\line(-1,-2){5}}
\end{picture}
\end{center}
\caption{Graphs of order $s\cdot (2k+1)$ for $k=5$ and $s\in \{ 1,2,5\}$. 
Considering the vertices of degree $2$, 
it follows easily that the encircled vertices form minimum $k$-component dominating sets.}\label{fig5}
\end{figure}
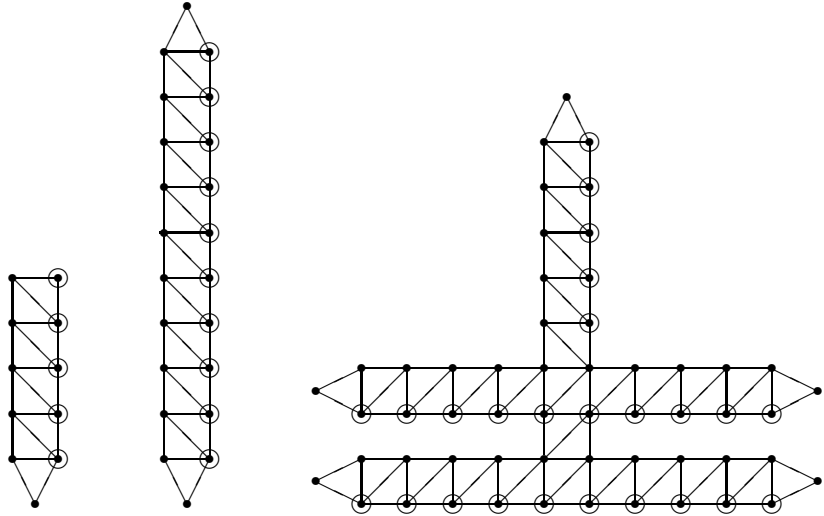

\noindent Let $k$ and $n$ be positive integers with $n\geq 2k+1$.
Lemma \ref{lemma3} and Lemma \ref{lemma5}(iii) imply that 
$$
\gamma_k(n)=
\left\{
\begin{array}{ll}
\lfloor\frac{kn}{2k+1}\rfloor, & \mbox{ if }n\in [2k+1,4k+3]\mbox{, and}\\
\lceil\frac{kn}{2k+1}\rceil, & \mbox{ if $n$ is an even number in }[4k+4,4k^2-2k].
\end{array}
\right.$$
Figure \ref{fig5} illustrates how to construct maximal outerplanar graphs $G$ 
of arbitrary order $n$ with $n\mod (2k+1)=0$ 
that satisfy $\gamma_k(G)=\lfloor\frac{kn}{2k+1}\rfloor$.

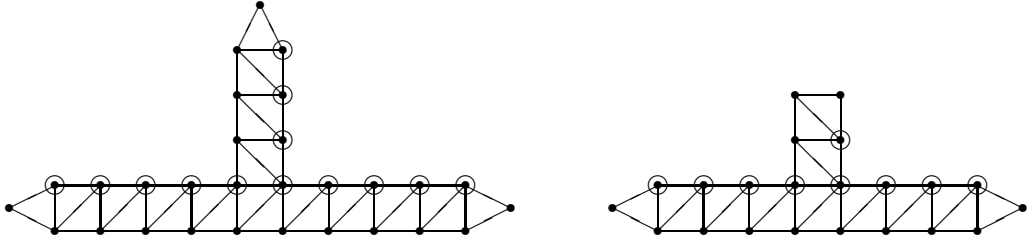
\begin{figure}[H]
\begin{center}
\unitlength 0.6mm 
\linethickness{0.4pt}
\ifx\plotpoint\undefined\newsavebox{\plotpoint}\fi 
\begin{picture}(111,51)(0,0)
\put(10,0){\circle*{2}}
\put(20,0){\circle*{2}}
\put(30,0){\circle*{2}}
\put(40,0){\circle*{2}}
\put(50,0){\circle*{2}}
\put(60,0){\circle*{2}}
\put(70,0){\circle*{2}}
\put(80,0){\circle*{2}}
\put(90,0){\circle*{2}}
\put(100,0){\circle*{2}}
\put(60,40){\circle*{2}}
\put(10,10){\circle*{2}}
\put(20,10){\circle*{2}}
\put(30,10){\circle*{2}}
\put(40,10){\circle*{2}}
\put(50,10){\circle*{2}}
\put(70,10){\circle*{2}}
\put(80,10){\circle*{2}}
\put(90,10){\circle*{2}}
\put(100,10){\circle*{2}}
\put(50,40){\circle*{2}}
\put(10,0){\line(0,1){10}}
\put(20,0){\line(0,1){10}}
\put(30,0){\line(0,1){10}}
\put(40,0){\line(0,1){10}}
\put(50,0){\line(0,1){10}}
\put(60,0){\line(0,1){10}}
\put(70,0){\line(0,1){10}}
\put(80,0){\line(0,1){10}}
\put(90,0){\line(0,1){10}}
\put(100,0){\line(0,1){10}}
\put(60,40){\line(-1,0){10}}
\put(10,10){\line(1,0){10}}
\put(20,10){\line(1,0){10}}
\put(30,10){\line(1,0){10}}
\put(40,10){\line(1,0){10}}
\put(50,10){\line(1,0){10}}
\put(70,10){\line(1,0){10}}
\put(50,20){\line(0,1){10}}
\put(80,10){\line(1,0){10}}
\put(90,10){\line(1,0){10}}
\put(50,30){\line(0,1){10}}
\put(20,10){\line(-1,-1){10}}
\put(30,10){\line(-1,-1){10}}
\put(40,10){\line(-1,-1){10}}
\put(50,10){\line(-1,-1){10}}
\put(60,10){\line(-1,-1){10}}
\put(70,10){\line(-1,-1){10}}
\put(80,10){\line(-1,-1){10}}
\put(50,30){\line(1,-1){10}}
\put(90,10){\line(-1,-1){10}}
\put(100,10){\line(-1,-1){10}}
\put(50,40){\line(1,-1){10}}
\put(10,0){\line(1,0){8}}
\put(20,0){\line(1,0){8}}
\put(30,0){\line(1,0){8}}
\put(40,0){\line(1,0){8}}
\put(50,0){\line(1,0){8}}
\put(60,0){\line(1,0){8}}
\put(70,0){\line(1,0){8}}
\put(60,20){\line(0,1){8}}
\put(80,0){\line(1,0){8}}
\put(90,0){\line(1,0){8}}
\put(60,30){\line(0,1){8}}
\put(18,0){\line(1,0){2}}
\put(28,0){\line(1,0){2}}
\put(38,0){\line(1,0){2}}
\put(48,0){\line(1,0){2}}
\put(58,0){\line(1,0){2}}
\put(60,18){\line(0,1){2}}
\put(68,0){\line(1,0){2}}
\put(78,0){\line(1,0){2}}
\put(88,0){\line(1,0){2}}
\put(98,0){\line(1,0){2}}
\put(60,38){\line(0,1){2}}
\put(10,10){\circle{4}}
\put(20,10){\circle{4}}
\put(30,10){\circle{4}}
\put(40,10){\circle{4}}
\put(50,10){\circle{4}}
\put(60,10){\circle{4}}
\put(60,20){\circle{4}}
\put(70,10){\circle{4}}
\put(80,10){\circle{4}}
\put(90,10){\circle{4}}
\put(100,10){\circle{4}}
\put(60,40){\circle{4}}
\put(60,20){\line(0,-1){11}}
\put(60,10){\line(-1,1){10}}
\put(50,20){\line(0,-1){10}}
\put(60,20){\circle*{2}}
\put(50,20){\circle*{2}}
\put(50,20){\line(1,0){10}}
\put(60,30){\circle*{2}}
\put(50,30){\circle*{2}}
\put(60,27){\line(0,1){3}}
\put(60,30){\line(-1,0){10}}
\put(60,30){\circle{4}}
\put(60,10){\line(1,0){10}}
\put(100,10){\line(2,-1){10}}
\put(110,5){\line(-2,-1){10}}
\put(60,40){\line(-1,2){5}}
\put(55,50){\line(-1,-2){5}}
\put(10,10){\line(-2,-1){10}}
\put(0,5){\line(2,-1){10}}
\put(55,50){\circle*{2}}
\put(110,5){\circle*{2}}
\put(0,5){\circle*{2}}
\put(60,10){\circle*{2}}
\end{picture}\hspace{1cm}
\linethickness{0.4pt}
\ifx\plotpoint\undefined\newsavebox{\plotpoint}\fi 
\begin{picture}(91,31)(0,0)
\put(10,0){\circle*{2}}
\put(20,0){\circle*{2}}
\put(30,0){\circle*{2}}
\put(40,0){\circle*{2}}
\put(50,0){\circle*{2}}
\put(60,0){\circle*{2}}
\put(70,0){\circle*{2}}
\put(80,0){\circle*{2}}
\put(10,10){\circle*{2}}
\put(20,10){\circle*{2}}
\put(30,10){\circle*{2}}
\put(40,10){\circle*{2}}
\put(50,10){\circle*{2}}
\put(60,10){\circle*{2}}
\put(70,10){\circle*{2}}
\put(80,10){\circle*{2}}
\put(10,0){\line(0,1){10}}
\put(20,0){\line(0,1){10}}
\put(30,0){\line(0,1){10}}
\put(40,0){\line(0,1){10}}
\put(50,0){\line(0,1){10}}
\put(60,0){\line(0,1){10}}
\put(70,0){\line(0,1){10}}
\put(80,0){\line(0,1){10}}
\put(10,10){\line(1,0){10}}
\put(20,10){\line(1,0){10}}
\put(30,10){\line(1,0){10}}
\put(40,10){\line(1,0){10}}
\put(50,10){\line(1,0){10}}
\put(40,20){\line(0,1){10}}
\put(60,10){\line(1,0){10}}
\put(70,10){\line(1,0){10}}
\put(20,10){\line(-1,-1){10}}
\put(30,10){\line(-1,-1){10}}
\put(40,10){\line(-1,-1){10}}
\put(50,10){\line(-1,-1){10}}
\put(60,10){\line(-1,-1){10}}
\put(40,30){\line(1,-1){10}}
\put(70,10){\line(-1,-1){10}}
\put(80,10){\line(-1,-1){10}}
\put(10,0){\line(1,0){8}}
\put(20,0){\line(1,0){8}}
\put(30,0){\line(1,0){8}}
\put(40,0){\line(1,0){8}}
\put(50,0){\line(1,0){8}}
\put(50,20){\line(0,1){8}}
\put(60,0){\line(1,0){8}}
\put(70,0){\line(1,0){8}}
\put(18,0){\line(1,0){2}}
\put(28,0){\line(1,0){2}}
\put(38,0){\line(1,0){2}}
\put(48,0){\line(1,0){2}}
\put(50,18){\line(0,1){2}}
\put(58,0){\line(1,0){2}}
\put(68,0){\line(1,0){2}}
\put(78,0){\line(1,0){2}}
\put(10,10){\circle{4}}
\put(20,10){\circle{4}}
\put(30,10){\circle{4}}
\put(40,10){\circle{4}}
\put(50,10){\circle{4}}
\put(50,20){\circle{4}}
\put(60,10){\circle{4}}
\put(70,10){\circle{4}}
\put(80,10){\circle{4}}
\put(50,20){\line(0,-1){11}}
\put(50,10){\line(-1,1){10}}
\put(40,20){\line(0,-1){10}}
\put(50,20){\circle*{2}}
\put(40,20){\circle*{2}}
\put(40,20){\line(1,0){10}}
\put(50,30){\circle*{2}}
\put(40,30){\circle*{2}}
\put(50,27){\line(0,1){3}}
\put(50,30){\line(-1,0){10}}
\put(10,10){\line(-2,-1){10}}
\put(0,5){\line(2,-1){10}}
\put(80,10){\line(2,-1){10}}
\put(90,5){\line(-2,-1){10}}
\put(0,5){\circle*{2}}
\put(90,5){\circle*{2}}
\end{picture}
\end{center}
\caption{On the left, a graph $G$ of order $n=s\cdot (2k+1)+(2t-1)$ for $k=5$, $s=2$, and $t=4$.
On the right, a graph $G$ of order $n=s'\cdot (2k'+1)+2t'$ for $k'=4$, $s'=2$, and $t'=2$.
Again, the encircled vertices form minimum $k$-component dominating sets.}\label{fig6}
\end{figure}

\noindent If $n\in \{ s\cdot (2k+1)+(2t-1), s\cdot (2k+1)+2t\}$ for positive integers $s$ and $t$ with $t\in [k]$, 
then $\lfloor\frac{kn}{2k+1}\rfloor=sk+t-1$.
The graphs in Figure \ref{fig6} illustrate how to construct extremal maximal outerplanar graphs of these orders.

\noindent Altogether, it follows that
$$
\gamma_k(n)=
\left\{
\begin{array}{ll}
\lceil\frac{kn}{2k+1}\rceil, & \mbox{ if $n$ is an even number in }[4k+4,4k^2-2k]\mbox{ and}\\
\lfloor\frac{kn}{2k+1}\rfloor, & \mbox{ otherwise}.
\end{array}
\right.$$

\medskip

\noindent We introduced the parameter $\gamma_k$ 
with the intention to obtain a common generalization of two separate results;
one concerning the domination number and one concerning the total domination number.
It seems interesting to unify/generalize further pairs of results about these parameters in this way.

\medskip

\noindent {\bf Acknowledgment}  
J.D. Alvarado and S. Dantas were partially supported by FAPERJ, CNPq, and CAPES.
D. Rautenbach was partially supported by CAPES.

\end{document}